\documentclass[reprint]{imsart}

\RequirePackage[OT1]{fontenc}
\RequirePackage{amsthm,amsmath,amssymb,amsfonts}
\RequirePackage[colorlinks,citecolor=blue,urlcolor=blue]{hyperref}

\arxiv{1512.07522}

\startlocaldefs
\numberwithin{equation}{section}

\theoremstyle{plain}
\newtheorem{theorem}{Theorem}[section]
\newtheorem{lemma}[theorem]{Lemma}
\newtheorem{remark}[theorem]{Remark}
\newtheorem{proposition}[theorem]{Proposition}
\newtheorem{definition}[theorem]{Definition}
\newtheorem{corollary}[theorem]{Corollary}
\newtheorem{example}[theorem]{Example}
\newtheorem{fig}[theorem]{Figure}

\usepackage{amsbsy}
\usepackage[T1]{fontenc}
\usepackage[utf8]{inputenc}
\usepackage{pgfplots}
\usepackage{wrapfig}
\usepackage{latexsym}
\usepackage{mathrsfs}
\usepackage{graphicx}
\usepackage{multicol}
\usepackage{multirow}
\usepackage{fancybox}
\usepackage{caption}
\usepackage{subcaption}
\usepackage{verbatim}

\usepackage{color,soul}
\usepackage[numbers]{natbib}
\usepackage{relsize}
\usepackage{lmodern}
\usepackage{url}
\usepackage{siunitx}
\usepackage{MnSymbol}
\usepackage[english]{babel}
\usepackage{tikz}
\usetikzlibrary{arrows}

\usepackage{algorithm}
\usepackage[noend]{algpseudocode}

\usepackage{caption} 

\newcommand{\bthe}{\begin{theorem}}
\newcommand{\ethe}{\end{theorem}}

\newcommand{\ben}{\begin{enumerate}}
\newcommand{\een}{\end{enumerate}}

\newcommand{\bit}{\begin{itemize}}
\newcommand{\eit}{\end{itemize}}

\newcommand{\beq}{\begin{equation}}
\newcommand{\eeq}{\end{equation}}

\newcommand{\ble}{\begin{lemma}}
\newcommand{\ele}{\end{lemma}}

\newcommand{\bde}{\begin{definition}\rm}
\newcommand{\ede}{\halmos\end{definition}}

\newcommand{\bco}{\begin{corollary}}
\newcommand{\eco}{\end{corollary}}

\newcommand{\bpr}{\begin{proposition}}
\newcommand{\epr}{\end{proposition}}

\newcommand{\brem}{\begin{remark}\rm}
\newcommand{\erem}{\halmos\end{remark}}

\newcommand{\bproof}{\begin{proof}}
\newcommand{\eproof}{\end{proof}}

\newcommand{\bexam}{\begin{example}\rm}
\newcommand{\eexam}{\halmos\end{example}}

\newcommand{\bfi}{\begin{fig}}
\newcommand{\efi}{\end{fig}}

\newcommand{\btab}{\begin{tab}}
\newcommand{\etab}{\end{tab}}

\newcommand{\beao}{\begin{eqnarray*}}
\newcommand{\eeao}{\end{eqnarray*}\noindent}

\newcommand{\beam}{\begin{eqnarray}}
\newcommand{\eeam}{\end{eqnarray}\noindent}

\newcommand{\ovr}{\begin{array}}
\newcommand{\earr}{\end{array}}

\newcommand{\bdis}{\begin{displaymath}}
\newcommand{\edis}{\end{displaymath}\noindent}
\newcommand{\MD}{{minimum {\rm ML} }}

\def\N{{\mathbb N}}

\def\R{{\mathbb R}}

\def\calb{{\mathcal{B}}}

\newcommand{\bs}{\boldsymbol}

\newcommand{\bfx}{\mathbf{X}}

\newcommand{\bfy}{\mathbf{Y}}

\def\1{\mathbf{1}}
\def\bone{\mathbf{1}}

\newcommand{\DAG}{{\rm DAG}}
\newcommand{\ML}{{\rm ML}}
\newcommand{\SEM}{{\rm SEM}}
\newcommand{\mSEM}{recursive {\rm ML}  model}

\newcommand{\CM}{{{\rm ML} coefficient matrix}}

\newcommand{\sgn}{{\rm sgn}}

\newcommand{\an}{{\rm an}}
\newcommand{\pa}{{\rm pa}}
\newcommand{\nd}{{\rm nd}}
\newcommand{\des}{{\rm de}}
\newcommand{\An}{{\rm An}}
\newcommand{\Pa}{{\rm Pa}}

\newcommand{\Des}{{\rm De}}
\newcommand{\tr}{{\rm tr}}
\newcommand{\tc}{{\rm tc}}

\newcommand{\low}{{\rm  low}}
\newcommand{\high}{{\rm  high}}

\newcommand{\nmw}{{\rm  nmw}}

\newcommand{\ov}{\overline}

\newcommand{\halmos}{\quad\hfill\mbox{$\Box$}}  

\def\D{\mathcal{D}}

\definecolor{plum}{cmyk}{0.50,1,0,0}
\definecolor{TealBlue}{cmyk}{0.86,0,0.34,0.02}
\definecolor{OliveGreen}{cmyk}{0.64,0,0.95,0.40}

\allowdisplaybreaks[4]
\endlocaldefs

\begin{document}

\begin{frontmatter}

\title{Max-linear models \\ on directed acyclic graphs}
\runtitle{Max-linear models on directed acyclic graphs}

\begin{aug}
  \author{\fnms{Nadine}  \snm{Gissibl}\corref{}\thanksref{t2}
  \ead[label=e1]{n.gissibl@tum.de}}
  \and
  \author{\fnms{Claudia} \snm{Kl\"uppelberg}\corref{}\thanksref{t2}
  \ead[label=e2]{cklu@tum.de} 
  \ead[label=u1,url]{http://www.statistics.tum.de}
 \ead[label=u1,url]{http://www.statistics.tum.de}}

  \thankstext{t2}{
{\sc  Address: } Center for Mathematical Sciences,  Boltzmannstrasse 3, 85748 Garching, Germany, 
         \printead{e1,e2,u1}}
         
\runauthor{Nadine Gissibl and Claudia Kl\"uppelberg}
\affiliation{Technische Universit\"at M\"unchen}
\address{Technical University of Munich 
  }
\end{aug}

\begin{abstract}
We consider a new recursive structural equation model where all  variables can be written as max-linear function of their parental node variables and independent noise variables.  
The model is max-linear in terms of the noise variables,  and its causal structure is represented by a directed acyclic graph.  
We detail the relation between the weights of the recursive structural equation model and the coefficients in its max-linear representation. 
In particular, we characterize all max-linear models which are generated by a recursive structural equation  model, and show that its max-linear coefficient matrix is the solution of a fixed point equation. 
We also find a unique minimum directed acyclic graph representing the recursive structural equations of the variables. 
The model structure introduces a natural order between the node variables and the max-linear coefficients. 
This yields representations of the vector components, which are based on a minimum number of node and noise variables.
\end{abstract}



\end{frontmatter}

\section{Introduction}\label{s1}

Graphical models are a popular tool to analyze and visualize the  conditional independence properties between random variables   (see e.g. \citet{KF} and \citet{Lauritzen1996}). 
Each node  in the graph indicates a random variable, and the absence of an edge between two nodes  represents a conditional independence property between the corresponding variables.  
We focus on directed graphical models, also called \textit{Bayesian networks}, where edge orientations come along with an intuitive causal interpretation.   
The conditional independence among the random variables, which is induced by a directed acylic graph (\DAG), can be explored using the (directed) Markov property:  each variable is conditionally independent of its non-descendants  (excluding the parents) given its parents  (cf. \cite{Lauritzen1996}, Chapter~3.2). 

Despite many areas of applications for  directed graphical models,  ranging from artificial intelligence, decision support systems, and engineering to genetics, geology, medicine, and finance (see e.g. \citet{Pourret2008}), graphical modelling of random vectors has mainly been limited to discrete and  Gaussian  distributions; see e.g. \cite{KF,Lauritzen1996}.  
In the context of risk assessment, risk exposures are usually modelled by continuous variables, however, the assumption of Gaussianity leads invariably to severe underestimation of large risks and therefore to unsuitable models.
  
Recursive structural equation models (recursive SEMs) offer a possibility to construct directed graphical models; cf. \citet{Bollen,Pearl2009} and \citet{SGS}. For a given \DAG\ $\D=(V,E)$ with nodes $V=\{1,\ldots,d\}$ and edges $E=\{(k,i) : i\in V\mbox{ and } k\in\pa(i)\}$ define
  \begin{align}\label{sem}
  X_i=f_i(\bfx_{{{\pa}}(i)},Z_i),\quad i=1, \dots, d,
  \end{align}
  where $\pa(i)$ denotes the parents of node $i$ in $\D$ and $f_i$ is a real-valued measurable function;
   $Z_1,\ldots,Z_d$ are independent noise variables. 
 Thus, a recursive SEM is specified by an underlying causal structure in terms of a \DAG\ $\D$, the functions $f_i$, and the distributions of $Z_i$  for $i=1,\ldots,d$.  
In this setting, the  distribution of $\bfx$ is uniquely defined by the  distributions of the noise variables and, denoting by ${\nd}(i)$ the non-descendants of node $i$,  
\begin{align}\label{mar}
X_i \upmodels \bfx_{{\nd}(i)\setminus {{\pa}}(i)} \mid\bfx_{{{\pa}}(i)}, \quad i=1, \dots, d;
\end{align}
i.e., the distribution of $\bfx$ is Markov relative to $\D$ (see Theorem~1.4.1 and the related discussion in \citet{Pearl2009}).  
Recently, recursive linear SEMs 
and generalisations in a Gaussian setting have received particular attention;  
see \citet{buhlmann2014cam,ernest2016causal} and references therein). 

Our focus is not on sums but on maxima,  where natural candidates for the noise distributions are the extreme value distributions or distributions in their domain of attraction (see e.g. \citet{Resnick1987,Resnick2007}).
We introduce a recursive \SEM, which is to the best of our knowledge new. 
Define a \textit{recursive max-linear} (\rm ML) \textit{model}  $\bfx=(X_1,\ldots,X_d)$ {\em on a} \DAG\  $\D$ by 
\begin{align}\label{ml-sem}
X_i:=\bigvee\limits_{k \in {{\pa}}(i)} c_{ki} X_k \vee c_{ii} Z_i,\quad i=1, \dots, d,
\end{align}
with independent non-negative random variables $Z_1, \dots, Z_d$
and positive weights $c_{ki}$ for all $i\in V$ and $k\in {\pa}(i)\cup\{i\}$. 

The new model is motivated by applications to risk analysis, where extreme risks play an essential role and may propagate through a network.  
In such a risk setting it is natural to require the noise variables to have positive infinite support $\R_+=[0,\infty)$. 
Moreover, we may think of the edge weights in \eqref{ml-sem} as relative quantities so that a risk may originate with certain proportions in its different ancestors.

In this paper we investigate  structural properties  as well as graph properties of a \mSEM\ $\bfx$ on a  \DAG\ $\D$.  
We will show that $\bfx$  is 
a \textit{max-linear {\em(ML)} model} (for background on \ML\ models in the context of extreme value theory see e.g. \citet{DHF}, Chapter~6) in the sense that
\begin{align}\label{maxlin}
X_i=\bigvee_{j=1}^d b_{ji} Z_j, \quad i=1,\ldots,d,
\end{align} 
with $Z_1,\ldots,Z_d$ as in \eqref{ml-sem},
and $B=(b_{ij})_{d\times d}$ is a matrix with non-negative entries.  
We call $B$ \textit{max-linear {\em(ML)} coefficient matrix} of $\bfx$ and its entries \textit{max-linear {\em(ML)} coefficients}.  

The \ML\ coefficients of $\bfx$ can be determined by a path analysis of $\D$. 
Throughout we  write  $k\to i$,  if there is  an edge from $k$ to $i$ in $\D$.   
We assign a weight to every path $p=[ j=k_0 \rightarrow k_1\rightarrow \dots \rightarrow  k_n=i]$, which is the product of the edge weights  along $p$ multiplied by the weight of the noise variable $Z_j$   (a concept, which goes back to \citet{wright}):
\begin{align}\label{bs}
d_{ji}(p) & = 
c_{k_0,k_0}  c_{k_0,k_1}  \dots c_{k_{n-2}, k_{n-1}} c_{k_{n-1},{k_n}}
 =   c_{k_0,k_0}\prod_{l=0}^{n-1} c_{k_l,k_{l+1}}.
\end{align} 
We will show that the \ML\ coefficients are given for $i\in V$ by  
\begin{align}\label{bs-max}
b_{ji}=\bigvee\limits_{p \in P_{ji}} d_{ji}(p)\mbox{ for $j \in {\an}(i)$};\, b_{ii}= c_{ii};\,\mbox{$b_{ji}=0$ for $j\in {V\setminus (\an(i)\cup\{i\}})$},
\end{align}  
where $P_{ji}$ is the set of paths from $j$ to $i$ and $\an(i)$  the ancestors of $i$. 


The computation in \eqref{bs-max} can be viewed as the algebraic path problem over the max-times semiring $(\R_+,\vee,\cdot)$ (see e.g. \citet{mahr1980birds} and \citet{rote1985systolic}). We present this problem in matrix form, using the matrix product over this semiring.
{We apply this concept in the two different situations, where the \DAG\ $\D$ of the model is given, and we test, if a given \ML\ coefficient matrix is consistent with $\D$,
but also later on, when we check, if a given matrix defines a recursive \SEM\ on some unspecified \DAG.}

From \eqref{bs-max} it is clear that not all paths are needed for representing $\bfx$ as ML model \eqref{maxlin}.
This perception leads to a complexity reduction of the model in different ways and in different situations.   
For every specific component $X_i$ of $\bfx$ only those paths with terminal node $i$, which carry the maximum weight, are relevant for its representation \eqref{maxlin}, and we call them \textit{max-weighted paths}.
All other paths can be disposed of without changing this representation.  
It is even sufficient to consider one max-weighted path in $\D$ from every ancestor of $i$ to $i$. 
Consequently, $X_i$ can be represented as component of a \mSEM\ on a polytree with node set $\An(i)$ 
and  with the same weights and noise variables as in the original representation \eqref{ml-sem}.

However, in general none of these individual polytrees represents all components of $\bfx$ in the sense of \eqref{ml-sem} simultaneously. 
Still there may be subgraphs of $\D$ and weights such that all components of $\bfx$ have representation \eqref{ml-sem}, and we present all such possible subgraphs and weights. 
In particular, we characterize the smallest subgraph $\D^B$ of this kind, which we call \textit{minimum max-linear {\em (ML)} \DAG\ of $\bfx$}, and point out its prominent role.

We are also interested in all  \DAG s, which represent $\bfx$ as a \mSEM, and show how the corresponding weights in representation \eqref{ml-sem} can be identified from the \CM\ of $\bfx$. 
In this context, we also give necessary and sufficient conditions on a matrix to be the \CM\ of any recursive \ML\ model. 

It is a simple but important observation that there is a natural order between the components of  $\bfx$;  from \eqref{ml-sem} we see immediately that $X_i\ge c_{ki} X_k$ holds for all $i\in V$ and $k\in\pa(i)$.
For every component  of $\bfx$ and some $U\subseteq V$, we find lower and upper bounds in terms of $\bfx_U:=(X_l,l\in U)$.  
Often we do not need all components of $\bfx_U$ to compute the best bounds of $X_i$ in terms of components of $\bfx_U$. 
If $i\in U$, then an upper and lower bound is given by $X_i$ itself;
otherwise, for a lower bound, we only need to consider a component $X_j$ of $\bfx_U$ if $j\in\an(i)$, but no max-weighted path from $j$ to $i$ passes through some node in $U\setminus \{j\}$.  
A similar result and  concept applies for the upper bound of $X_i$. 
Thus, the max-weighted paths also lead in this context indirectly to a complexity reduction. 
We will also use the max-weighted ancestors of $i$ in $U$ to obtain a minimal representation of $X_i$ in terms of $\bfx_U$ and noise variables.

 
Our paper is organized as follows. In Section~\ref{s2} we discuss the max-linearity of a \mSEM\ $\bfx$ and express  its \CM\  in terms of a weighted adjacency matrix of a corresponding \DAG. 
Section~\ref{s3} introduces the important notion of  a max-weighted path and studies its consequences for the \ML\ coefficients. 
In Section~\ref{s4}  we give necessary and sufficient conditions for a \ML\ model being a \mSEM\ on a given \DAG. Section~\ref{s5} is devoted to the minimum \ML\ \DAG\ of $\bfx$ as the \DAG\ with the minimum number of edges within the class of all \DAG s representing $\bfx$ in the sense of \eqref{ml-sem}. 
In Section~\ref{s6}, given a  set of node variables, we investigate which information can be drawn for the other components of $\bfx$. 
This results in lower and upper bounds for the components. 
Finally, we derive a minimal representation for the components of $\bfx$ as max-linear functions of a subset of node variables and certain noise variables. 
  
\smallbreak
 We use the following notation throughout. 
For a node $i \in V$, the sets ${\an}(i)$, ${{\pa}}(i)$, and  ${\des}(i)$ contain the {\em ancestors, parents}, and {\em descendants} of $i$ in $\D$.  
Furthermore, we use the notation $\An(i)={\an}(i)\cup\{i\}$,  ${\Pa}(i)={{\pa}}(i)\cup\{i\}$, and  ${\Des}(i)={\des}(i)\cup\{i\}$. 
We write $U\subseteq V$ for a non-empty subset $U$ of nodes, $\bfx_U=(X_l, l\in U)$, and $U^c=V\setminus U$. 
All our vectors are row vectors. 
We also extend the previous notation in a natural way by writing ${\an}(U)=\bigcup_{i\in U}\an(i)$,  $\An(U)={\an}(U)\cup U$, and so on. 
For a matrix $B$ with non-negative entries we write $\sgn(B)$ for the matrix with entries equal to 1, if the corresponding component in $B$ is positive and 0 else.
We denote by $\bone_U$ the indicator function of $U$, and set $\bone_{\emptyset}\equiv 0$.  
In general, we consider statements for $i \in \emptyset$ as invalid. 
Moreover, for arbitrary (possibly random)  $a_i\in\R_+$,  we set $\bigvee_{i\in\emptyset}  a_i =0$  and  $\bigwedge_{i\in\emptyset}  a_i =\infty$.

\section{Max-linearity of a recursive max-linear model\label{s2}}

For a recursive \ML\ model $\bfx$ on a \DAG\ $D=(V,E)$, given by \eqref{ml-sem}, we derive its max-linear representation \eqref{maxlin}.
We start with our leading example, the diamond-shaped \DAG\  depicted  below. 

\bexam\label{ex:maxlin}[Max-linear representation of a \mSEM]\\
Consider a \mSEM\ $\bfx=(X_1,X_2,X_3,X_4)$ with \DAG
\begin{align*}
\D=(V,E)=(\{1,2,3,4\}, \{(1,2),(1,3),(2,4),(3,4)\})
\end{align*}
and weights $c_{ki}$ for $i\in V$ and $k\in\Pa(i)$.
We obtain for the random variables $X_1,$ $X_2$, $X_3$, and $X_4$:
\begin{multicols}{2}
\vspace*{-1cm}
\begin{minipage}{0.5\textwidth}
\begin{center}
\begin{align*}
X_1 &= c_{11} Z_1\\
X_2 &= c_{12} X_1\vee c_{22} Z_2 = c_{12} c_{11}Z_1\vee  c_{22} Z_2\\
X_3 &= c_{13} X_1\vee c_{33} Z_3 = c_{13}c_{11} Z_1 \vee  c_{33} Z_3\\
X_4 &= c_{24} X_2 \vee c_{34} X_3\vee  c_{44} Z_4 \\
&=c_{24} (c_{12}c_{11} Z_1\vee c_{22} Z_2) \vee c_{34}(c_{13}c_{11}Z_1 \vee c_{33} Z_3)\vee  c_{44} Z_4\\
&=  ( c_{24}c_{12}c_{11}  \vee c_{34}c_{13}c_{11}) Z_1 \vee c_{24}  c_{22}Z_2 \vee c_{34} c_{33} Z_3\vee c_{44} Z_4.
\end{align*}
\end{center}
\end{minipage}
\begin{minipage}{0.5\textwidth}
\flushright
\begin{tikzpicture}[->,every node/.style={circle,draw},line width=1pt, node distance=1.8cm,minimum size=0.9cm]
 \node (1) [outer sep=1mm] {$1$};
  \node (2) [below left of=1,outer sep=1mm] {$2$};
    \foreach \from/\to in {1/2}
  \draw (\from) -- (\to);   
  \node (3) [below right of=1,outer sep=1mm] {$3$};
    \foreach \from/\to in {1/3}
  \draw (\from) -- (\to);   
  \node (4) [below right of=2,outer sep=1mm] {$4$}; 
    \foreach \from/\to in {2/4,3/4}
  \draw (\from) -- (\to);   
\end{tikzpicture}
\end{minipage}
\end{multicols}
\noindent
Thus $\bfx$ satisfies \eqref{maxlin} with \CM
\begin{align*}
B = 
\begin{bmatrix}
c_{11} &  c_{11}c_{12}    & c_{11} c_{13}   &  c_{11}c_{12} c_{24}  \vee c_{11}c_{13}c_{34}   \\
0  & c_{22} &   0    &  c_{22}c_{24} \\
0    &    0   &   c_{33} &  c_{33} c_{34}  \\
0 & 0 & 0 &   c_{44}
\end{bmatrix};
\end{align*}
i.e., the \ML\ coefficients satisfy \eqref{bs-max}. 
Moreover, $B$ is an upper triangular matrix, since $\D$ is well-ordered (cf. Remark~\ref{maxlinear_struc}(ii)).
\eexam

The following result shows that such a representation can be obtained in general: every component of a  \mSEM\ has a max-linear representation in terms of its ancestral noise variables  and an independent one. 
It also provides a general method to calculate the \ML\ coefficients by a path analysis as described in \eqref{bs} and \eqref{bs-max}. 

\bthe\label{cor:MaxLinRep}
Let $\bfx$ be a \mSEM\ with \DAG\ $\D$, and let $B=(b_{ij})_{d\times d}$ be the  matrix with entries as defined in \eqref{bs-max}.   
Then  
\begin{align}\label{ml-noise}
X_i 
=\bigvee\limits_{j \in \An(i)} b_{ji}Z_j, \quad i=1,\ldots,d;
\end{align}
i.e., $B$ is the \ML\ coefficient matrix of $\bfx$.
\ethe

\bproof
We know that every \DAG\  may be well-ordered (see Remark~\ref{maxlinear_struc}(ii)).
Hence, without loss of generality we assume throughout this proof that $\D$ is well-ordered. 
We prove the identity \eqref{ml-noise} by induction on the number of nodes of  $\D$. 
For $d=1$ we have by \eqref{ml-sem}
\begin{align*}
X_1= c_{11} Z_1=b_{11} Z_1,
\end{align*}
where the last equality holds by \eqref{bs-max}.
Suppose  that \eqref{ml-noise} holds for a \mSEM\ $\bfx$ of dimension $d$; i.e., 
\begin{align*}
X_k =  \bigvee_{j \in \An(k)} b_{jk} Z_{j}=  \bigvee_{j \in {\an}(k)} \bigvee\limits_{p \in P_{jk}}d_{jk}(p) Z_{j} \vee c_{kk} Z_k,\quad k=1,\ldots, d.
\end{align*}
Now consider a $(d+1)$-variate \mSEM, and note that for every $i\in \{1,\ldots,d\}$ we have $(d+1)\in V\setminus{{\pa}}(i)$, since $\D$ is well-ordered. 
Thus,  in order to verify \eqref{ml-noise} for the nodes $i\in\{1,\dots,d\}$, it suffices to consider the subgraph $\D[\{1,\dots,d\}]= (\{1,\dots,d\}, E\cap (\{1,\dots,d\}\times \{1,\dots,d\} ))$. 
 Due to the induction hypothesis, \eqref{ml-noise} holds for $\D[\{1,\dots,d\}]$ and, hence, also for $\D$. 
So we can use  this hypothesis and \eqref{lem11} to obtain  
\begin{align*} 
 X_{d+1}&= \bigvee\limits_{k\in {{\pa}}(d+1)} c_{k,d+1} X_k \vee c_{d+1,d+1} Z_{d+1}\\
  &= \bigvee\limits_{k\in {{\pa}}(d+1)}   \bigvee\limits_{j \in {\an}(k)} \bigvee\limits_{p \in P_{jk}} c_{k,d+1} d_{jk}(p) Z_j \vee \bigvee\limits_{k\in {{\pa}}(d+1)}  c_{k,d+1} c_{kk} Z_k   \vee c_{d+1,d+1} Z_{d+1}\\
  &= \bigvee_{j \in {\an}(d+1)} \big(  \bigvee_{k \in   {\des}(j)\cap {{\pa}}(d+1)}  \bigvee_{p \in P_{jk}} c_{k,d+1} d_{jk}(p) \vee 
  \bigvee_{k \in {{\pa}}(d+1)\cap \{j\}} c_{k,d+1}c_{kk} \big) Z_j  \\
& \quad\quad \vee c_{d+1,d+1} Z_{d+1}.
  \end{align*}
Observe that every path from some $j$ to $d+1$ is of the form $p=[j\to\ldots\to k\to d+1]$ for some $k\in{\des}(j)\cap {{\pa}}(d+1)$, or an edge $j\to d+1$ corresponding to $j\in\pa(d+1)$. 
From \eqref{bs}, the path $p$  has weight $d_{j,d+1}(p)= d_{jk}(p)c_{k,d+1}$, and the edge $j\to d+1$ has weight  $d_{j,d+1}([j\to d+1])=c_{jj}c_{j,d+1}$.
This yields 
 \begin{align*}  
X_{d+1}    &= \bigvee_{j \in {\an}(d+1)} \bigvee_{p \in P_{j,d+1}} d_{j,d+1}(p) Z_j    \vee c_{d+1,d+1} Z_{d+1}
=\bigvee_{j \in \An(d+1)} b_{j,d+1}Z_j,
\end{align*}
where we have used that  $b_{j,d+1}=\bigvee\limits_{p \in P_{j,d+1}}d_{j,d+1}(p)$ for $j \in {\an}(d+1)$   and $b_{d+1,d+1}= c_{d+1,d+1}$. 
\eproof

By \eqref{bs-max} the \ML\ coefficient $b_{ji}$ of $\bfx$ is different from zero if and only if $j\in \An(i)$.
This information is contained in the  {\em reachability matrix}  $R=(r_{ij})_{d \times d}$ of $\D$, which has entries
\begin{align*}
r_{ji}:=\begin{cases}
1, & \text{if there is a path from $j$ to $i$, or if $j=i$},\\
0, & \text{otherwise}. 
\end{cases}
\end{align*}
If the $ji$-th entry of $R$ is equal to one, then $i$ {\em is reachable from} $j$.

\brem\label{maxlinear_struc} 
Let $\D$ be a \DAG\ with reachability matrix $R$.
\begin{enumerate} 
\item[(i)]
The \CM\ $B$ is a weighted reachability matrix of  $\D$; i.e., $R=\sgn(B)$.
\item[(ii)]
Every \DAG\ $\D$ can be \textit{well-ordered}, which means that the set
$V=\{{1}, \dots,d\}$ of nodes is linearly ordered in a way compatible with $\D$  such that $k\in{{\pa}}(i)$ implies $k<i$  (see e.g. Appendix A of \citet{Diestel:2010}).
If $\D$ is well-ordered, then $B$ and $R$ are upper triangular matrices.
\end{enumerate}
\vspace*{-0.5em}
\erem

Finding the \CM\ $B$ from $\D$ and the weights in \eqref{ml-sem} by  a path analysis as described in \eqref{bs} and \eqref{bs-max} would be very inefficient. 
We may, however, compute $B$ by means of a specific matrix multiplication.   

For two non-negative matrices $F$ and $G$, where the number of columns in $F$ is equal to the number of rows in $G$, we define the product $\odot: \mathbb{R}_+^{m\times n}\times \mathbb{R}_+^{n\times p} \rightarrow \mathbb{R}_+^{m\times p}$ by
\begin{align}\label{odot} 
(F=(f_{ij})_{m\times n},G=(g_{ij})_{n\times p}) \mapsto F\odot G :=\Big(\bigvee\limits_{k=1}^n f_{ik}g_{kj}\Big)_{m\times p}. 
\end{align}
The triple  $(\R_+,\vee,\cdot)$, which is called max-times or subtropical algebra, is an idempotent semiring with $0$ as 0-element and  $1$ as 1-element. 
The operation $\odot$ is therefore a matrix product over a semiring. 
 Such semirings are fundamental in tropical geometry; for an introduction see \citet{MS}.  
The matrix product $\odot$ is associative: for $F\in \R_+^{ m\times n } $, $G\in  \R_+^{ n\times p} $, and  $H\in  \R_+^{ p\times q} $, $F\odot (G\odot H)=(F \odot G) \odot H$,  and we have $(F\odot G)^\top = G^\top\odot F^\top$. 
Denoting by $\calb$ all $d\times d$ matrices with non-negative entries and by $\vee$ the componentwise maximum between two matrices, $(\calb,\vee,\odot)$ is also a semiring with the null matrix as 0-element and the identity matrix  ${\rm id}_{d\times d}$  as 1-element.  
This semiring is, however, not commutative, since $\odot$ is in general not. 
Consistent with a matrix product,  we define powers recursively: $A^{\odot 0} := {\rm id}_{d\times d}$ and  
$A^{\odot n} := A^{\odot (n-1)}\odot A$ for $A\in\calb$ and $n\in\N$.  

The matrix product $\odot$ allows us to present the problem of characterising representation \eqref{ml-noise} from \eqref{ml-sem} in terms of $B$, involving the weighted adjacency matrix $(c_{ij} {\bone}_{\pa(j)}(i))_{d\times d}$ of $\D$. 


\bthe\label{the:odot}
Let $\bfx$ be a \mSEM\   with  \DAG\ $\D$ and weights  $c_{ki}$ for $i\in V$ and $k\in\Pa(i)$ as in \eqref{ml-sem}.
Furthermore, define the matrices  
\begin{align*}
A:= {\rm diag}(c_{11},\ldots,c_{dd}),\, A_0:= 
\big(c_{ij} \1_{\pa(j)}(i)\big)_{d\times d},\,\mbox{and}\,\, A_1:=
\big(c_{ii} c_{ij} \1_{\pa(j)}(i) \big )_{d\times d}.
\end{align*} 
Then the \CM\  $B$ of $\bfx$ {from Theorem~\ref{cor:MaxLinRep}} has representation
\begin{align*}
B = A \quad \text{for $d=1$} \quad \text{and} \quad B =  A \boldsymbol{\vee} \bigvee_{k=0}^{d-2} \big(A_1 \odot A_0^{\odot k}\big) \quad\text{for $d\ge 2$}.
\end{align*}
\ethe

\bproof
For $d=1$ we know from  \eqref{bs-max} that $b_{11}=c_{11}$. Hence, $B=A$. 
Now assume that $d\ge 2$.  
 First we show that, if $\D$ has a path of length $n$  (a path consisting of $n$ edges) from node $j$ to node $i$, then the $ji$-th entry of the matrix  $A_1\odot A_0^{\odot(n-1)}$ is equal to the maximum weight of all paths of lengths $n$ from $j$ to $i$, otherwise it is zero.  The proof is by induction on $n$. 

An edge $j\to i$, which  is the only path of length $n=1$, has the weight $d_{ji}([j\to i]) = c_{jj} c_{ji}$. 
Since the $ji$-th entry of the matrix $A_1\odot A_0^{\odot 0}=A_1\odot  {\rm id}_{d\times d}=A_1$  is given by $c_{jj}c_{ji}{\bs 1}_{\pa(i)}(j)$, the statement is true for $n=1$. 

Denote by $a_{n,ji}$ and $a_{n+1,ji}$ the $ji$-th entry of $A_1\odot A_0^{\odot (n-1)}$ and $A_1\odot A_0^{\odot n}$, respectively.    
As $A_1\odot A_0^{\odot n}=(A_1\odot A_0^{\odot (n-1)})\odot A_0$, the   $ji$-th entry   of $A_1\odot A_0^{\odot n}$ is given by     $a_{n+1,ji}=\bigvee_{k=1}^{d} a_{n,jk}a_{0,ki}=\bigvee_{k=1}^{d} a_{n,jk}c_{ki} \1_{\pa(i)}(k)$.  We obtain from the induction hypothesis and \eqref{bs} that $a_{n,jk}a_{0,ki}$  is zero, if $\D$ does not contain a path of length $n$ from $j$ to $k$ or  the edge $k\to i$; otherwise it is  equal to the maximum weight of all paths which consist of a path of length $n$ from $j$ to $k$ and the edge $k\to i$. 
Since every path of length $n+1$ from $j$ to $i$ is of this form for some $k\in V$, the $ji$-th entry of $A_1\odot A_0^{\odot n}$ is indeed equal to the maximum weight of all paths of length $n+1$ from $j$ to $i$ if there exists such a path, otherwise it is zero. 

Finally, recall from \eqref{bs-max} that for $i\in V$ and $j\in\an(i)$ the \ML\ coefficient $b_{ji}$ is equal to the maximum weight of all paths from $j$ to $i$, and note that due to acyclicity,  a path  in $\D$ is at most of length $d-1$. 
Thus, if $j\in\an(i)$ then the $ji$-th entry of $\bigvee_{k=0}^{d-2}A_1\odot A_0^{\odot k}$  is equal to  $b_{ji}$, otherwise it is zero. 
Since by \eqref{bs-max}, $b_{ii}=c_{ii}$ and $b_{ji}=0$ for $j\in V\setminus\An(i)$, the \CM\  $B$  is given by
\begin{align*}
B = A \vee A_1 \vee \big(A_1 \odot A_0 \big) \vee \big(A_1 \odot A_0^{\odot 2}\big)\vee \dots\vee  \big(A_1\odot A_0^{\odot (d-2)} \big).
\end{align*}
\eproof

The following has been shown in the proof of Theorem~\ref{the:odot}.

\bco\label{cor2.5}
If $\D$ has a path  of length $n$ from $j$ to $i$, the  $ji$-th entry of the matrix $A_{1}\odot A_{0}^{\odot(n-1)}$ is equal to the maximum weight of  all paths of length $n$ from $j$ to $i$, otherwise the entry is zero. 
\eco

Summarizing the noise variables of $\bfx$  into the  vector $\mathbf{Z}=(Z_1,\ldots,Z_d)$, the representation \eqref{ml-noise}  of $\bfx$ can be written by means of the product $\odot$ as
\begin{align*} 
\bfx=\mathbf{Z} \odot B= \big(\bigvee_{j=1}^d b_{ji} Z_j, i=1,\ldots,d\big) = \big(\bigvee_{j\in\An(i)} b_{ji} Z_j, i=1,\ldots,d\big).
\end{align*}
Consequently, the definition of the matrix product~$\odot$ modifies and extends the definition given in \citet[Section~2.1, Eq. (2)]{Wang2011}.

\section{Max-weighted paths and submodels \label{s3}}

Given a \mSEM\ $\bfx$ with \DAG\ $\D=(V,E)$, weights $c_{ki}$ for $i\in V$ and  $k\in\Pa(i)$, and \CM\ $B=(b_{ij})_{d\times d}$, we investigate the paths of $\D$, their particular weights, relations between the \ML\ coefficients, and an induced subgraph structure.

From \eqref{bs-max} and \eqref{ml-noise} we know that a path $p$ from $j$ to $i$, whose weight $d_{ji}(p)$ is strictly smaller than $b_{ji}$ does not have any influence on the distribution of $\bfx$. 
This fact suggests the following definition. 

\bde\label{mwp}
Let $\bfx$ be a recursive \ML\ model with \DAG\ $\D=(V,E)$, \ML\ coefficient matrix $B$, and path weights as in \eqref{bs}.
We call a path $p$ from $j$ to $i$ a  {\em max-weighted path (in $\D$)} if  $b_{ji}=d_{ji}(p)$.  
\ede

A prominent example, where all paths are max-weighted, is the following.

\bexam\label{polytree}[Polytree]
A {\em polytree} is a \DAG\ whose underlying undirected graph has no cycles;   
polytrees have at most one path between any pair of nodes. 
Thus, assuming  that $\bfx$ is a \mSEM\ on a polytree, all paths  must be max-weighted.
\eexam

The next example emphasizes the importance and consequences of max-weighted paths, which we will investigate in more detail in the next sections.

\bexam\label{Breduce}[Max-weighted path, graph reduction]\\ 
Consider   a \mSEM\ $\bfx=(X_1,X_2,X_3)$ with \DAG
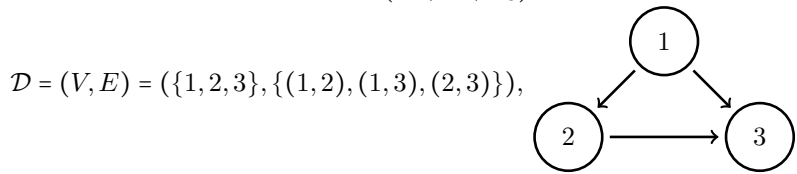
 \begin{multicols}{2}
  \begin{minipage}{0.6\textwidth}
\begin{center}
  \begin{align*}
 \D=(V,E)= (\{1,2,3\}, \{(1,2),(1,3),(2,3)\}),
  \end{align*}
  \end{center}
  \end{minipage}
\begin{minipage}{0.4\textwidth} 
\flushright
\vspace*{-1em}
\begin{tikzpicture}[->,every node/.style={circle,draw},line width=1pt, node distance=1.8cm,minimum size=0.9cm]
 \node (1) [outer sep=1mm]  {$1$};
  \node (2) [below left of=1,outer sep=1mm] {$2$};
  \node (3) [below right of=1,outer sep=1mm] {$3$};
    \foreach \from/\to in {1/2,1/3,2/3}
  \draw (\from) -- (\to);   
\end{tikzpicture}
\end{minipage}
\end{multicols}
\noindent weights $c_{11}, c_{12}, c_{13}, c_{22}, c_{23}, c_{33}$, and  \CM\ $B$. 
We distinguish between two situations: \\[1mm]
(1)\, If  $c_{13}> c_{12} c_{23}$, then the edge $1\to 3$ is the unique max-weighted path from $1$ to $3$. \\[1mm]
(2)\,  If, however,  $c_{13}\le  c_{12} c_{23}$, then  $b_{13}=c_{11}c_{12}c_{23}=\frac{b_{12} b_{23}}{b_{22}}$ and  the path $[1\to 2\to 3]$ is max-weighted. We obtain in this case
\begin{align*}
 X_3 &= b_{13}Z_1\vee b_{23}Z_2\vee b_{33}Z_3=\frac{b_{23}}{b_{22}} (b_{12}Z_1\vee b_{22}Z_2) \vee b_{33}Z_3  = c_{23}X_2\vee b_{33}Z_3. 
\end{align*}
Thus,  $\bfx$ is also a \mSEM\ on the  \DAG\ 
\begin{align*}
\D^B:=(\{1,2,3\}, \{(1,2),(2,3)\}).
\end{align*}
Here $\D^B$ is the  \DAG\ with minimum number of edges such that $\sgn(B)$ is its reachability matrix.
\eexam 

We present some immediate consequences of the path weights in \eqref{bs} and the definition of max-weighted paths. 

\brem\label{rem:mwp}
\begin{enumerate} 
\item[(i)] 
If there is only one path between two nodes, it is max-weighted. 
\item[(ii)]  
Every subpath of a max-weighted path is also max-weighted. 
\item[(iii)] 
Every path, which results from a max-weighted path by replacing a subpath with another max-weighted subpath, is also max-weighted. 
\end{enumerate}
\vspace*{-0.5cm}
\erem

To find for some $i\in V$ and $j\in \an(i)$ the \ML\ coefficient $b_{ji}$ it suffices to know the weight of $Z_j$ and the edge weights along one arbitrary max-weighted path from $j$ to $i$, since every max-weighted path from $j$ to $i$ has  the same weight. 
This allows us to represent every component $X_i$ of $\bfx$ as component of a \mSEM\ on a subgraph of $\D$.  
For this purpose, we introduce the following definition. 

\bde \label{def:inducedmod}
Let $\overline D=(\ov V, \ov E)$ be a subgraph of $\D$, and denote by $\ov \pa(i)$ the  parents of node $i$ in $\ov\D$.  
Define 
\begin{align*}
Y_i:=\bigvee_{k\in\ov\pa(i)}c_{ki}Y_k\vee c_{ii}Z_i, \quad i\in \ov V,
\end{align*}
with the same weights and noise variables as for $X_i$ in representation \eqref{ml-sem}. We   call the resulting \mSEM\ $\bfy=(Y_l,l\in\ov V)$ \textit{recursive \ML\ submodel of $\bfx$ induced by $\ov\D$}. 
\ede

We summarize some immediate properties of $\bfy$.

\brem\label{rem:indsubMod} 
Let $i\in V$ with ancestors $\an(i)$ in $\D$.
Denote by $\ov B=(\ov b_{ij})_{\vert \ov V\vert\times \vert \ov V\vert}$ the \CM\ of  $\bfy$. 
\begin{enumerate}
\item[(i)] Every path  in $\ov\D$ has the same weight \eqref{bs} as in $\D$. 
\item[(ii)] A path of $\ov\D$, which is in $\D$ a max-weighted path, is also in $\ov\D$  max-weighted. 
\item[(iii)] For $j\in\an(i)$, $\ov\D$ has one in $\D$ max-weighted path from $j$ to $i$ if and only if $\ov b_{ji}=b_{ji}$. 
\item[(iv)] $\ov \D$ has at least one in $\D$ max-weighted path  from every $j\in\an(i)$ to $i$ if and only if $X_i=Y_i$. 
\end{enumerate}
\vspace*{-0.8em}
\erem

By Remark~\ref{rem:mwp}(ii),  for every $i\in V$, there exists a polytree $\D_i$ of $\D$ with node set $\An(i)$, which has exactly one in $\D$ max-weighted path from every ancestor of $i$ to $i$. 
There may even exist several such polytrees (cf. Example~\ref{ex:ml} below).  
We learn  from the construction of $\D_i$ and Remark~\ref{rem:mwp}(ii) that indeed every path of $\D_i$ is in $\D$ max-weighted. 
Therefore, some component $X_j $ of $\bfx$ coincides  by Remark~\ref{rem:indsubMod}(iv) with the corresponding one of the recursive \ML\ submodel of $\bfx$ induced by $\D_i$ if and only if $\D_i$ has at least one path from every ancestor of $j$ in $\D$ to $j$. 
By construction of $\D_i$ this property holds obviously for $X_i$. We summarize this result as follows.

\bpr \label{le:polytrees2}
Let $\bfx$ be a recursive \ML\ model with \DAG\ $\D$ and \ML\ coefficient matrix $B$.
For some $i\in V$ and $\An(i)$ in $\D$ let $\D_i$ be a polytree with node set $\An(i)$ such that $\D_i$ has one in $\D$ max-weighted path from every $j\in\an(i)$ to $i$. 
Let $\bfy_i=(Y_l, l\in \An(i))$ be the recursive \ML\ submodel of $\bfx$ induced by  $\D_i$. 
Then for all $j\in \An(i)$, which have the same ancestors in $\D_i$ and $\D$, we have $X_j=Y_j$.
\epr

We discuss the \mSEM\ from  Example~\ref{ex:maxlin}  in the context of Definition~\ref{mwp} and Proposition~\ref{le:polytrees2}. 

\bexam\label{ex:ml}[Continuation of Example~\ref{ex:maxlin}: max-weighted paths,  polytrees, conditional independence]\\
We identify all max-weighted paths ending in node $4$.  
By Remark~\ref{rem:mwp}(i), the paths $[2\to 4]$ and $[3\to 4]$  are max-weighted. 
For the weights of the paths from node $1$ to $4$ we have three situations:
\begin{align*}
c_{11}c_{12}c_{24}=c_{11}c_{13}c_{34},\quad c_{11}c_{12}c_{24}>c_{11}c_{13}c_{34},\quad  \text{and} \quad c_{11}c_{12}c_{24}<c_{11}c_{13}c_{34}.
\end{align*}
In the first situation, both paths from $1$ to $4$, $[1\to 2 \to 4]$ and $[1\to 3 \to 4]$, are max-weighted. 
Thus, there are two different polytrees having one in $\D$ max-weighted  path from every ancestor of $4$ to $4$, namely,
\begin{align*}
\D_{4,1} &=(\{1,2,3,4\},\{(1,2),(2,4),(3,4)\}) \quad \text{and} \quad \\
\D_{4,2} &=(\{1,2,3,4\},\{(1,3),(2,4),(3,4)\}).
\end{align*}
In the second situation,  the path $[1\to 2 \to 4]$ is the unique max-weighted path from $1$ to $4$ and, hence, $\D_{4,1}$ is the unique polytree as in Proposition~\ref{le:polytrees2} for node 4. 
The third case is symmetric to the second, such that  $\D_{4,2}$ is also such a unique polytree. 

Now let $\bfy_1=(Y_{1,1},Y_{1,2},Y_{1,3},Y_{1,4})$ and $\bfy_2=(Y_{2,1},Y_{2,2},Y_{2,3},Y_{2,4})$ be the recursive \ML\ submodels of $\bfx$ induced by  $\D_{4,1}$ and $\D_{4,2}$. If the path $[1\to 2 \to 4]$  is max-weighted, we have by Proposition~\ref{le:polytrees2} that
\begin{align}\label{ex:ml1}
Y_{1,1}=X_1, \quad Y_ {1,2}=X_2, \quad  \text{and} \quad  Y_ {1,4}=X_4; 
\end{align}
if $[1\to 3 \to 4]$  is max-weighted, then
\begin{align}\label{ex:ml2}
Y_{2,1}=X_1, \quad  Y_ {2,3}=X_3, \quad \text{and} \quad Y_{2,4}=X_4. 
\end{align}
We know  that the distributions of $\bfx$, $\bfy_1$, and $\bfy_2$ are Markov relative to $\D$, $\D_{4,1}$, and  $\D_{4,2}$, respectively.   
For a \DAG, the local Markov property as specified in \eqref{mar}, is {by Proposition~4 of  \citet{Lauritzen1990}} equivalent to the global Markov property (for a definition see Corollary~3.23 of  \cite{Lauritzen1996}).
Using this property we find 
 \begin{align*}
 Y_{1,1}\upmodels Y_{1,4} \mid  Y_{1,2} \quad \text{and}  \quad Y_{2,1}\upmodels Y_{2,4} \mid  Y_{2,3} . 
 \end{align*}
Thus, if the path $[1\to 2 \to 4]$  is in $\D$ max-weighted, we have  by \eqref{ex:ml1}  that $X_{1}\upmodels X_4 \mid  X_{2}$. Accordingly, if $[1\to 3 \to 4]$ is max-weighted,  $X_{1}\upmodels X_4 \mid  X_{3}$ holds by \eqref{ex:ml2}. 
Since the only conditional independence property encoded in $\D$ by the (global) Markov property is  $X_1 \upmodels X_4\mid X_2, X_3$, we can identify additional conditional independence properties of $\bfx$ from the polytrees in Proposition~\ref{le:polytrees2}. 
\eexam

\brem 
(i)\, Assume the situation of Proposition~\ref{le:polytrees2}. 
Let $V_i$ be the set of all  nodes in $\An(i)$, which have the same ancestors in $\D$ and $\D_i$. 
Since the distributions of $\bfx$ and $\bfy$ are Markov relative to $\D$ and $\D_i$, respectively, conditional independence properties of $\bfx$ are encoded in $\D$ and of $\bfy$ in $\D_i$.  
By  Proposition~\ref{le:polytrees2}, the conditional independence properties between subvectors of  $\bfy_{V_i}=(Y_l,l\in V_i)$, which we can read off from $\D_i$,  hold also between the corresponding subvectors of $\bfx$. 
Since missing edges correspond to conditional independence properties, 
and $\D_i$ is a subgraph of $\D$, we can often identify additional  conditional independence properties of $\bfx$  from $\D_i$.  \\[1mm]
(ii)\, From (i) or Example~\ref{ex:ml} we learn that a \mSEM\  with \DAG\ $\D$ is in general not \textit{faithful}; i.e., not all conditional independence properties are  encoded in $\D$ by the (global) Markov property.  
\erem

 As can be seen from  Example~\ref{ex:ml}, any reduction of a \mSEM\ depends on the existence of max-weighted paths that pass through some specific node.  
The following result  shows how we can obtain this information from its \CM.

\bthe\label{paththrough}
Let  $B$  be the \ML\ coefficient matrix of a \mSEM\  on the \DAG\ $\D$.  
Let further $U\subseteq V$, $i \in V$ and $j\in \an(i)$, and recall from Remark~\ref{maxlinear_struc}(i) that $b_{ji}>0$.   
\begin{enumerate}
\item[(a)] There is a max-weighted path from $j$ to $i$, which passes through some node in $U$ 
if and only if 
\begin{align}
b_{ji} &=\bigvee_{k \in {\Des}(j)\cap U\cap \An(i)} \frac{b_{jk}b_{ki}}{b_{kk}}.\label{paththrough2} 
\end{align}
\item[(b)]  No max-weighted path from $j$ to $i$ passes through some node in $U$ if and only if
\begin{align}
b_{ji}&> \bigvee_{k \in {\Des}(j)\cap U \cap {\An}(i)} \frac{b_{jk}b_{ki}}{b_{kk}}.   \label{paththrough3} 
\end{align}
This holds also for $U=\emptyset$.
\end{enumerate}
\ethe

\bproof
First assume that $\Des(j)\cap U\cap\An(i)=\emptyset$. 
Thus no path, hence also no max-weighted path, from $j$ to $i$ passes through some node in $U$, and it suffices to verify (b).  
Since the right-hand side of \eqref{paththrough3} 
is zero if and only if $\Des(j)\cap U\cap\An(i)=\emptyset$, and the \ML\ coefficient $b_{ji}$ is positive, (b) is proven for this case. 

Now assume that  $\Des(j)\cap U\cap\An(i)=\{k\}$, which implies that there is a path from $j$ to $i$ passing through $k\in U$.   
If $k=i$ or $k=j$, there is obviously a max-weighted path from $j$ to $i$ passing through $i$ or $j$ and \eqref{paththrough2} is always valid. 

Next assume that $k\in V\setminus\{i,j\}$ and that $p_1$ as well as $p_2$ are max-weighted paths from $j$ to $k$ and from $k$ to $i$. 
Denote by $p$ the path from $j$ to $i$ consisting of the subpaths $p_1$ and $p_2$.
By \eqref{bs} and the definition of a max-weighted path we obtain
\begin{align*}
d_{ji}(p)=\frac{1}{c_{kk}}d_{jk}(p_1) d_{ki}(p_2) =\frac{b_{jk}b_{ki}}{b_{kk}}. 
\end{align*}
Since  $p$ is max-weighted if and only if $b_{ji}=d_{ji}(p)$, and this is not the case if and only if  $b_{ji}>d_{ji}(p)$, we have shown (a) and (b) for the situation of $\Des(j)\cap U\cap\An(i)=\{k\}$.  
In particular, it follows  that $b_{ji}\ge \frac{b_{jk}b_{ki}}{b_{kk}}$ for all $k\in \Des(j)\cap U\cap \An(i)$. 

Assume now  that  $\Des(j)\cap U\cap\An(i)$ contains more than one element, and that a max-weighted path from $j$ to $i$ passes through some node $k\in U$. We know from above that this is equivalent to
\begin{align*}
b_{ji}=\frac{b_{jk}b_{ki}}{b_{kk}} \quad  \text{and} \quad b_{ji} \ge \frac{b_{jl}b_{li}}{b_{ll}} \,\,\, \text{for all $l\in  \big(\Des(j)\cap U\cap\An(i)\big)\setminus\{k\}$}, 
\end{align*}
which is again equivalent to \eqref{paththrough2}. Similarly, we obtain (b). 
\eproof

\brem \label{rem2:paththrough}
Recall the matrix product~$\odot$ from \eqref{odot}.  We obtain from $R=\sgn(B)$ (Remark~\ref{maxlinear_struc}(i)) that for $i,j\in V$
\begin{align*}
\bigvee_{k \in {\Des}(j)\cap U \cap {\An}(i)} \frac{b_{jk}b_{ki}}{b_{kk}}=\bigvee_{k=1}^d \frac{b_{jk}b_{ki}}{b_{kk}} \1_{U}(k)
=:\bigvee_{k=1}^d b_{jk} b_{U,ki}
\end{align*}
 is the $ji$-th entry of the matrix $B\odot B_U$ with  $B_U=(b_{U,ij})_{d\times d}$. 
Thus, we may decide whether there is a max-weighted path between two nodes that passes through some node in $U$ by comparing the entries of the matrices $B$ and $B\odot B_U$.  
Such use of the matrix product~$\odot$ can  be made  at various points throughout the paper, for instance in Remark~\ref{rem:DB}(i), Theorem~\ref{lem:DB},~and~Lemma~\ref{poss1}(b). 
\erem

From Theorem~\ref{paththrough}, recalling from Remark~\ref{maxlinear_struc}(a) that $R=\sgn(B)$, we obtain an important property of the \ML\ coefficients. 

\bco \label{rem:paththrough}   
For all $i\in V$, $k\in \An(i)$, and $j\in \An(k)$, $b_{ji} \ge \frac{b_{jk}b_{ki}}{b_{kk}}>0$. 
Indeed, $b_{ji}\ge \frac{b_{jk}b_{ki}}{b_{kk}}$ holds for all $i,j,k\in V$.
\eco

We learn immediately from \eqref{ml-sem} that $c_{ki} X_k\le X_i$ for all  $i\in V$ and $k\in{\pa}(i)$.  
From Corollary~\ref{rem:paththrough} we find such inequalities also for components, whose nodes are not connected by an edge but by a path of arbitrary length.

\bco \label{le:IneqXi}
For all $i\in V$ and $j\in \An(i)$, we have $\frac{b_{ji}}{b_{jj}} X_j \le X_i$. 
\eco

\bproof
Note that $\An(j)\subseteq\An(i)$. 
Using the max-linear representation \eqref{ml-noise} of $X_i$ and $X_j$ as well as Corollary~\ref{rem:paththrough}, we obtain
\begin{align*}
X_i=\bigvee_{l\in \An(i)} b_{li}Z_l \ge \bigvee_{l\in \An(j)} b_{li}Z_l \ge \bigvee_{l\in \An(j)} \frac{b_{lj}b_{ji}}{b_{jj}}Z_j=  \frac{b_{ji}}{b_{jj}} \bigvee_{l\in \An(j)} b_{lj}Z_j=\frac{b_{ji}}{b_{jj}} X_j.
\end{align*}
\eproof

\section{ML coefficients leading to a recursive \ML\ model on a given DAG\label{s4}}

Recall the definition of a (general) \ML\ model given in \eqref{maxlin}.
From Theorem~\ref{cor:MaxLinRep} we know that every \mSEM\ is  max-linear.
In this section we provide necessary and sufficient conditions on a \ML\ model  to be  a \mSEM\ on a given \DAG\ $\D$. 

It can be shown that every \ML\ model, which is a recursive \SEM\ 
as given in \eqref{sem} with unspecified functions $f_1,\dots,f_d$, must be a \mSEM.
That a \mSEM\ on $\D$ is also a recursive \SEM\  follows immediately from its recursive definition. 
To summarize, a  \ML\ model can be represented as a recursive SEM \eqref{sem} with \DAG\ $\D$ if and only if it has a recursive ML representation \eqref{ml-sem} relative to the same \DAG\ $\D$.

We investigate below, when a \ML\ coefficient matrix $B$ as in \eqref{maxlin} is the \CM\ of a \mSEM\ with given \DAG\ $\D$. 
Motivated by Remark~\ref{maxlinear_struc}(i) in what follows we assume that $\sgn(B)$ is the reachability matrix $R$ of $\D$.
In our investigation the \DAG\ with the minimum number of edges, such that $R=\sgn(B)$, 
will play an important role. 
This  has already been indicated in  Example~\ref{Breduce}.

We give a general definition of the \DAG\ with minimum number of edges that represents the same reachability relation as a given \DAG.

\bde \label{Dtr}
Let $\D=(V,E)$ be a \DAG. The \DAG\ $\D^{\tr}=(V,E^{\tr})$ is the {\em transitive reduction} of $\D$ if the following holds:
\begin{enumerate}
\item[(a)] 
$\D^{\tr}$ has a path from node $j$ to node $i$ if and only if $\D$ has a path from $j$ to $i$, and
\item[(b)] there is no graph with less edges than $\D^{\tr}$ satisfying condition (a).
\end{enumerate}
\vspace*{-0.5cm}
\ede 

Since we work with finite \DAG s throughout, the transitive reduction is unique and is also a subgraph of the original \DAG.
The transitive reduction of a \DAG\ can be obtained by successively examining its edges, in any order, and deleting an edge $k\to i$, if the original \DAG\ contains a path from $k$ to $i$ which does not include this edge.
For these properties and further details see e.g. \citet{Aho1972}.
In what follows we need the notion of ${{\pa}}^{\tr}(i)$, the parents of $i$ in $\D^\tr$.

We present necessary and sufficient conditions on $B$ to be the \ML\ coefficient matrix of a recursive \ML\ model on $\D$. 

\bthe\label{fixedpointeq}
 Let $\D$ be a  \DAG\ with reachability matrix $R$ and $\bfx$ a \ML\ model as in \eqref{maxlin} with \CM\  $B$ such that ${\sgn}(B)=R$.
 Define 
  \begin{align*}
 A:={\rm diag}(b_{11},\ldots,b_{dd}) \quad\mbox{and}\quad
 A_0
 := \Big(\frac{b_{ij}}{b_{ii}}{\1}_{\pa(j)}(i) \Big)_{d \times d}.
 \end{align*}
  Then $\bfx$ is a \mSEM\ on $\D$ if and only if 
  {the following fixed point equation holds:}
\begin{align}\label{fpe}
B=A \vee B \odot A_0.
\end{align}
In this case, 
\begin{align*}
X_i=\bigvee_{k\in {{\pa}}(i)}\frac{b_{ki}}{b_{kk}}X_k \vee b_{ii}Z_i,\quad i=1,\dots,d.
\end{align*}
\ethe

\bproof 
First we investigate the fixed point equation \eqref{fpe} and compute the $ji$-th entry of $B\odot A_0$. 
By definition,  together with ${\sgn}(B)=R$,  it is equal to 
\begin{align*}
 \bigvee_{k=1}^d  \frac{b_{jk}b_{ki}}{b_{kk}}{\1}_{\pa(i)}(k) = \bigvee_{k\in {\Des}(j)\cap {{\pa}}(i)}  \frac{b_{jk}b_{ki}}{b_{kk}}.
\end{align*} 
We have ${\Des}(j)\cap {{\pa}}(i)=\emptyset$ for $j \in V\setminus\an(i)$ and ${\Des}(j)\cap {{\pa}}(i)={\des}(j)\cap {{\pa}}(i)$ for $j\in\an(i)\setminus \pa(i)$. Moreover, for $j \in {{\pa}}^{\tr}(i)$,  using that  ${\des}(j)\cap {{\pa}}(i)=\emptyset$, we obtain ${\Des}(j)\cap {{\pa}}(i)=\{i\}$. Thus, taking also the  matrix $A$ into account, \eqref{fpe} is equivalent  to 
\begin{align*}
b_{ji}=
\begin{cases}
0, & \text{if $j \in V\setminus\An(i)$,}\\
b_{ii}, & \text{if $j= i$,}\\
 \bigvee\limits_{k\in {\des}(j)\cap {\pa}(i)}  \dfrac{b_{jk}b_{ki}}{b_{kk}},\quad & \text{if $j \in {\an}(i)\setminus {\pa}(i)$,}\\
 b_{ji}\vee \bigvee\limits_{k\in {\des}(j)\cap {\pa}(i)}  \dfrac{b_{jk}b_{ki}}{b_{kk}},\quad & \text{if $j \in {\pa}(i)\setminus {\pa}^{\tr}(i)$},\\
 b_{ji}, & \text{if $j \in {\pa}^{\tr}(i)$}
\end{cases}
\end{align*}
for all $i,j\in V$. 
To summarize,  the fixed point equation \eqref{fpe} is   satisfied  if and only if for all $i\in V$ the following identities hold:
 \begin{align}
 b_{ji} &= \bigvee_{k \in   {\des}(j)\cap {\pa}(i)  }\frac{b_{jk}  b_{ki}}{b_{kk}}\quad\quad\quad \text{for all } j \in {\an}(i)\setminus  {\pa}(i) \label{Coeffeq1a}\\
b_{ji} &=  b_{ji} \vee \bigvee_{k \in   {\des}(j) \cap {\pa}(i) } \frac{b_{jk}   b_{ki}}{b_{kk}} \quad \text{for all } j \in {\pa}(i)\setminus {\pa}^{\tr}(i).\label{Coeffeq2a}
\end{align}
Thus it suffices to show that, under the conditions above,  $\bfx$ is a \mSEM\ on $\D$ if and only if \eqref{Coeffeq1a} and \eqref{Coeffeq2a} hold for all  $i\in V$.

First assume that $\bfx$ is a \mSEM\ on $\D$, and let $i\in V$ and $j\in\an(i)$. 
 Since  
 every path from $j$ to $i$  passes through at least one parent node of $i$, there must be a max-weighted path from $j$ to $i$ passing through some node in $\pa(i)$. 
Using \eqref{paththrough2} with $U=\pa(i)$ and noting that $j\in \Des(j)\cap U\cap\An(i)= \Des(j)\cap \pa(i)$, we find for $j\in\an(i)\setminus\pa(i)$ Eq. \eqref{Coeffeq1a} and for $j\in\pa(i)\setminus\pa^\tr(i)$ Eq. \eqref{Coeffeq2a}. 

For the converse statement, assume that \eqref{Coeffeq1a} and \eqref{Coeffeq2a} hold. For $j\in\pa^\tr(i)$  we have   $\des(j)\cap \pa(i)=\emptyset$, such that  the right-hand side of  \eqref{Coeffeq2a}  is equal to $b_{ji}$. Thus  \eqref{Coeffeq2a}  holds for all $j\in\pa(i)$.  Since $\sgn(B)=R$, we have $X_i=\bigvee_{j=1}^db_{ji}Z_j=\bigvee_{j\in\An(i)}b_{ji}Z_j$.  We split up the index set and use \eqref{Coeffeq1a} in the first place and \eqref{Coeffeq2a} for all $j\in\pa(i)$  in the second place to obtain
\begin{align*}
X_{i}
&= \bigvee_{j \in {\an}(i)\setminus {\pa}(i)} b_{ji}Z_j \vee \bigvee_{j \in  {\pa}(i)} b_{ji}Z_j  \vee b_{ii}Z_i\\
&=\bigvee_{j \in {\an}(i)\setminus {\pa}(i)}  \bigvee_{k \in   {  {\des}(j) \cap \pa}(i) }\frac{b_{jk}  b_{ki}}{b_{kk}} Z_j \vee \bigvee_{j \in  {\pa}(i)} b_{ji} Z_j \\
& \quad\quad\quad \vee \bigvee_{j \in  {\pa}(i)} \bigvee_{k \in {\des}(j) \cap \pa(i)  } \frac{b_{jk}   b_{ki}}{b_{kk}}Z_j  \vee b_{ii}Z_i\\
&=\bigvee_{j \in {\an}(i)}  \bigvee_{k \in   {  {\des}(j) \cap \pa}(i) }\frac{b_{jk}  b_{ki}}{b_{kk}} Z_j \vee \bigvee_{j \in  {\pa}(i)} b_{ji} Z_j  \vee b_{ii}Z_i.
\end{align*}
Interchanging the first two maximum operators by \eqref{lem11} yields 
\begin{align*}
X_{i}
&=\bigvee_{k\in{\pa}(i)} \bigvee_{j\in\an(k)} \frac{b_{jk}  b_{ki}}{b_{kk}} Z_j \vee \bigvee_{k \in  {\pa}(i)} b_{ki} Z_k  \vee b_{ii}Z_i\\
&= \bigvee_{k \in {\pa}(i) }  \frac{ b_{ki}}{b_{kk}}  \big(  \bigvee_{j \in {\an}(k)} b_{jk}Z_j 
\vee  b_{kk} Z_k \big) \vee b_{ii}Z_i \\
&= \bigvee_{k \in {\pa}(i) }  \frac{ b_{ki}}{b_{kk}}  X_k\vee b_{ii}Z_i.
\end{align*}
\eproof

In the proof of Theorem~\ref{fixedpointeq} we have shown that under the required conditions  the fixed point equation \eqref{fpe} holds if and only if  \eqref{Coeffeq1a} and \eqref{Coeffeq2a} hold. 
We summarize this  in part (a) of the following corollary.  
Part (b) has also been verified  in the proof of Theorem~\ref{fixedpointeq}.
The final statement is based on the fact that for $k\in\pa(i)$ we have $\des(k)\cap\pa(i)=\emptyset$ if and only if $k\in\pa^\tr(i)$.

\bco \label{equiDML}
(a) \,  Assume the situation of Theorem~\ref{fixedpointeq}.  \\
Then $\bfx$ is a \mSEM\  on $\D$ if {and only if} for every $i\in V$,
 \begin{align} 
 b_{ji} &= \bigvee_{k \in   {\des}(j)\cap {\pa}(i)  }\frac{b_{jk}  b_{ki}}{b_{kk}}\quad \text{for all } j \in {\an}(i)\setminus  {\pa}(i)\label{Coeffeq1}\\
b_{ji} &\geq  \bigvee_{k \in   {\des}(j) \cap {\pa}(i) } \frac{b_{jk}   b_{ki}}{b_{kk}} \quad \text{for all } j \in {\pa}(i)\setminus {\pa}^{\tr}(i) \label{Coeffeq2}.
\end{align}
(b) \, Let $\bfx$ be a recursive \ML\ model with \DAG\ $\D$ and \ML\ coefficient matrix $B$.
Then for every $i\in V$ and  $k\in\pa(i)$,
\begin{align*}
 b_{ki} \ge \bigvee_{l \in   {\des}(k)\cap {\pa}(i)  }\frac{b_{kl}  b_{li}}{b_{ll}}.
\end{align*}
Moreover, the right-hand side is equal to 0 if and only if $k\in\pa^\tr(i)$.
\eco

By \eqref{Coeffeq1} and \eqref{Coeffeq2} exactly those \ML\ coefficients $b_{ki}$, such that $k\to i$ is an edge in $\D^\tr$, do not have to meet any specific conditions apart from being positive. 

In summary, given a \DAG\ $\D$ with  node set $V=\{1,\ldots,d\}$, both Theorem~\ref{fixedpointeq}  and Corollary~\ref{equiDML}(a) characterize all \ML\ coefficient matrices of any \mSEM\ possible on  $\D$ as all non-negative $d\times d$  matrices that are weighted reachability matrices of $\D$ and satisfy  \eqref{fpe}, equivalently \eqref{Coeffeq1} or \eqref{Coeffeq2}. 
If we can verify these two properties for a non-negative $d\times d$  matrix $B$, then it is the \CM\ of a \mSEM\ on $\D$, and  for $i\in V$ weights in  its representation \eqref{ml-sem} are  given by $c_{ki}=\frac{b_{ki}}{b_{kk}}$ for $k\in \pa(i)$ and $c_{ii}=b_{ii}$.

\section{Graph reduction for a recursive max-linear model\label{s5}}

From Proposition~\ref{le:polytrees2} we know that every component $X_i$ of a recursive \ML\ model $\bfx$ with \DAG\ $\D=(V,E)$ satisfies \eqref{ml-sem}  on a subgraph of $\D$. 
These subgraphs, however, usually vary from one vector component to another. 
On the other hand,  we know from Example~\ref{Breduce} that the whole vector $\bfx$ may  also be a \mSEM\ on a subgraph of  $\D$.  
This raises the question of finding the smallest subgraph of $\D$ such that $\bfx$ is a \mSEM\ on this \DAG. 
We define and characterize this unique minimal \DAG\  before we  point out its prominent role in the class of all \DAG s representing $\bfx$ in the sense of \eqref{ml-sem}. 

\bde\label{DB}
Let $\bfx$ be a \mSEM\ with \DAG\ $\D=(V,E)$ and \CM\ $B$. 
We call the \DAG \ 
\begin{align} \label{DB1}
\D^B=(V, E^B):=\Big(V, \big\{(k,i) \in E:   b_{ki}> \bigvee_{l\in \des(k)\cap \pa(i)} \frac{b_{kl}b_{li}}{b_{ll}} \big\}\Big)
\end{align}
the {\em minimum max-linear {\em(ML)} \DAG\ of $\bfx$}. 
\ede

We summarize some properties of $\D^B$ as follows. 

\brem\label{rem:DB}
(i)\, By Theorem~\ref{paththrough}(b) the \MD \DAG\ $\D^B$ contains exactly those edges {$k\to i$} of  $\D$, where no max-weighted path from $k$ to $i$ passes through some node in $\pa(i)\setminus \{k\}$. 
This means that $\D^B$ has an edge $k\to i$ if and only if it is the only max-weighted path from $k$ to $i$ in $\D$. 
The \DAG\ $\D^B$ can be obtained from $\D$ by deleting an edge $k\to i$, if $\D$ contains a max-weighted path from $k$ to $i$, which does not include this edge. 
The algorithm is by comparison of the \ML\ coefficients:
for all $i\in V$ and $k\in \pa(i)\setminus \pa^{\tr}(i)$ remove the edge $k\to i$ from $\D$ if 
\begin{align*}
 b_{ki}=\bigvee_{l \in {\des}(k) \cap {\pa}(i)} \frac{b_{kl}   b_{li}}{b_{ll}}.
\end{align*}
Note the analogy to finding the transitive reduction $\D^{\tr}$ of $\D$ below Definition~\ref{Dtr}.\\[1mm]
(ii) The  \MD \DAG\ $\D^B=(V,E^B)$ is a subgraph of the original \DAG\ $\D=(V,E)$.  
Recall that the transitive reduction  $\D^\tr=(V,E^\tr)$ of $\D$ is also a subgraph of $\D$ and that every edge $k\to i$ in $\D^\tr$ is the only -- and hence also max-weighted -- path from $k$ to $i$ in $\D$.  
Thus, the transitive reduction $\D^\tr$ is also a subgraph of $\D^B$. In summary, we have $E^\tr \subseteq E^B \subseteq E$.
This implies that the \DAG s $\D^B$ and $\D$ have the same reachability matrix $\sgn(B)$.
\erem

The method described in Remark~\ref{rem:DB}(i)  determines $\D^B$ from  $\D$ and  $B$. 
Indeed, we can also identify  $\D^B$ directly from $B$ without knowing $\D$.

\bthe\label{lem:DB}
Let $\bfx$ be a \mSEM\  with \CM\ $B$. 
Then the  \MD \DAG\  of $\bfx$ can be represented as 
\begin{align}\label{eq:dagDB}
\D^B=\Big(V, \Big\{(k,i)\in V\times V: k\neq i \text{ and }   b_{ki}> \bigvee_{l=1\atop l\neq i,k}^d \frac{b_{kl}b_{li}}{b_{ll}}\Big\}\Big);
\end{align}
in particular, $\D^B$ is identifiable from $B$.
\ethe

\bproof 
Let $\D$ be a \DAG, which describes $\bfx$ in the sense of \eqref{ml-sem}.
{Since $R=\sgn(B)$ (Remark~\ref{maxlinear_struc}(i)) we have
\begin{align}\label{char:DB}
b_{ki}> \bigvee_{l=1\atop {l\neq i,k}}^d \frac{b_{kl}b_{li}}{b_{ll}}=\bigvee_{l\in \des(k)\cap \an(i)}  \frac{b_{kl}b_{li}}{b_{ll}}.
\end{align}
We show that the edge set in \eqref{eq:dagDB} coincides with $E^B$ as defined in \eqref{DB1}. 
Assume first that $(k,i)$ is contained in the edge set in \eqref{eq:dagDB}. 
Such a \DAG\  exist by the definition of a \mSEM.  }
Since  the right-hand side of \eqref{char:DB} is non-negative, we must have $b_{ki}>0$ and, hence, $k\in \an(i)$. 
By Theorem~\ref{paththrough}(b) no max-weighted path from $k$ to $i$ passes through some node in $V\setminus \{i,k\}$. 
Thus the edge $k\to i$ must be the only max-weighted path from $k$ to $i$ and, hence, by Remark~\ref{rem:DB}(i) it must be an edge  $E^B$ as in \eqref{DB1}. 

For the converse, let $(k,i)\in E^B$. 
Since by Remark~\ref{rem:DB}(i) this edge is the only max-weighted path from $k$ to $i$, there is no max-weighted path passing through some node in $V\setminus \{i,k\}$. 
This is by Theorem~\ref{paththrough}(b) equivalent to \eqref{char:DB} and $(k,i)$ belongs to the edge set in \eqref{eq:dagDB}.
\eproof

We  characterize all \DAG s and specify all weights such that  $\bfx$ satisfies \eqref{ml-sem}.
The \MD\DAG\ $\D^B$  of $\bfx$  is the smallest \DAG\ of this kind and  has unique weights in representation \eqref{ml-sem} in the sense that all irrelevant weights are set to zero. We can add edges into  $\D^B$ with weights $c_{ki}\in (0,\frac{b_{ki}}{b_{kk}}]$ representing $\bfx$ again in the sense of \eqref{ml-sem} as long as the graph represents the same reachability relation as $\D^B$. 
As a consequence, to find $B$ by a path analysis as described in \eqref{bs-max} it suffices to know $\D^B$  and the weights in representation \eqref{ml-sem} relative to $\D^B$. 


\bthe\label{le1:CMgraph} 
Let $\bfx$ be a \mSEM\ with   \CM\ $B$.  Let further $\D^B=(V,E^B)$ be the  \MD \DAG\  of $\bfx$ and $\pa^B(i)$ be the parents of node $i$ in  $\D^B$. 
\begin{enumerate}
\item[(a)] The  \MD \DAG\ $\D^B$ of $\bfx$ is  the \DAG\  with the minimum number of edges such that $\bfx$ satisfies \eqref{ml-sem}. 
The weights in \eqref{ml-sem} are uniquely given by $c_{ii}=b_{ii}$ and $c_{ki}=\frac{b_{ki}}{b_{kk}}$  for $i\in V$ and $k\in\pa^B(i)$. 
\item[(b)] Every  \DAG\  with node set $V$ that has at least the edges of $\D^B$ and the same reachability matrix as $\D^B$ represents $\bfx$ in the sense of \eqref{ml-sem} with weights 
 given for all $i\in V$ by
 \begin{align*}
 c_{ii}=b_{ii}; \, c_{ki}=\frac{b_{ki}}{b_{kk}}\text{ for $k\in  \pa^B(i)$}, \text{ and } c_{ki}\in \Big(0,\frac{b_{ki}}{b_{kk}}\Big]\text{ for $k\in \pa(i)\setminus \pa^B(i)$}.  
 \end{align*}
There are no further \DAG s and   weights   such that  $\bfx$ has representation \eqref{ml-sem}. 
\end{enumerate}
\ethe

\bproof 
(a)\, { Let $\D$ be a \DAG\ and $c_{ki}$ for $i\in V$ and $k\in\Pa(i)$ weights such that $\bfx$ has representation \eqref{ml-sem}. By Remark~\ref{rem:DB}(ii) $\D^B$ is a subgraph of $\D$. }

First we prove  that $\bfx$ is a \mSEM\ on $\D^B$  with weights $c_{ki}$ for $i\in V$ and $k\in\Pa^B(i)$ by showing that  all components of $\bfx$ coincide with those of the recursive \ML\ submodel of $\bfx$ induced by $\D^B$ (see Definition~\ref{def:inducedmod}).  
By  Remark~\ref{rem:indsubMod}(iv), it suffices to verify  for all $i\in V$ and $j\in \an(i)$ that $\D^B$ has  one in $\D$ max-weighted path from $j$ to $i$.  Among all  max-weighted paths from $j$ to $i$ in $\D$, let $p$ be one with maximal length, and
assume that  $p$ includes an edge, say $k\to l$, which is not contained in $\D^B$. 
The \DAG\ $\D$ has by Remark~\ref{rem:DB}(i), however,  a max-weighted path $p_1$ from $k$ to $l$, which does not include the edge $k \to l$. 
Note that $p_1$ consists of more edges than the path $[k\to l]$. Thus by replacing in $p$ the edge $k \to l$ by $p_1$ we obtain by Remark~\ref{rem:mwp}(iii) a max-weighted path from $j$ to $i$  consisting of more edges than $p$. Since this  a contradiction to the fact that  $p$ has maximal length among all max-weighted paths from $j$ to $i$,  $p$ must be in $\D^B$. 

Since every edge $k\to i$ in $\D^B$ is by Remark~\ref{rem:DB}(ii)  the only  max-weighted path from $k$ to $i$ in $\D$, we have by Definition~\ref{mwp} and  \eqref{bs} that $b_{ki}=c_{kk}c_{ki}=b_{kk}c_{ki}$, which implies $c_{ki}=\frac{b_{ki}}{b_{kk}}$, and  these  weights  are uniquely given. For the same reason  there cannot be  a \DAG\ such that $\bfx$ has representation \eqref{ml-sem} with less edges than  $\D^B$.  

(b)\,  From Remark \ref{le1:CMgraph}(ii) every \DAG\ $\D$ that represents  $\bfx$ in the sense of \eqref{ml-sem} must  have  the same reachability matrix as $\D^B$ and must contain at least  the edges of $\D^B$.  By  \eqref{bs} and \eqref{bs-max} the weights in representation \eqref{ml-sem} of $\bfx$ have to satisfy  $c_{ki}\le \frac{b_{ki}}{b_{kk}}$ for all $i\in V$ and $k\in\pa(i)$. 

It remains to show that  $\bfx$  satisfies \eqref{ml-sem} relative to  a \DAG\ $\D$ with the properties and weights $c_{ki}$ for $i\in V$ and $k\in \Pa(i)$ (the parents in $\D$)  as in the statement of (b). 
Note that  the \DAG\ $\D^B$ is a subgraph of $\D$ and both \DAG s have the same reachability relation.  
Since $\bfx$ is by part  (a) a \mSEM\ on $\D^B$,  we may   use Corollary~\ref{le:IneqXi} with the ancestors in  $\D^B$:  for every $i\in V$ and $k\in\pa(i)$, since $k$ is an ancestor of $i$ in $\D^B$ and $\frac{b_{ki}}{b_{kk}}\ge c_{ki}$, we have
\begin{align*}
X_i \ge \frac{b_{ki}}{b_{kk}}X_k \ge c_{ki}X_k.
\end{align*}
With this we obtain from representation \eqref{ml-sem} of $X_i$ relative to $\D^B$ that
\begin{align*}
X_i=\bigvee_{k\in\pa^B(i)} c_{ki}X_i\vee c_{ii}Z_i=\bigvee_{k\in\pa^B(i)}c_{ki} X_k\vee \bigvee_{k\in\pa(i)\setminus\pa^B(i)}c_{ki}X_k \vee c_{ii} Z_i,
\end{align*}
which is \eqref{ml-sem} relative to $\D$ with weights $c_{ki}$ for $i\in V$ and $k\in \Pa(i)$.
\eproof

As explained before Theorem~\ref{le1:CMgraph} we can add edges into $\D^B$, while keeping the same reachability relation and still having representation \eqref{ml-sem} for $\bfx$.
In what follows we will use  the \DAG\ with the maximum number of edges with these properties. 

\bde \label{Dtc}
Let $\D=(V,E)$ be a \DAG. The {\em transitive closure} $\D^{\tc}=(V,E^{\tc})$  of  $\D$ is the \DAG\ with edge $j\to i$  if and only if $\D$ has a path from  $j$ to  $i$. 
\ede 

The transitive reduction is essentially the inverse operation of the transitive closure:  for the transitive reduction one reduces the number of edges and for the transitive closure one adds edges, while maintaining the identical reachability relation. The transitive reduction of a \DAG\ $\D$ is a subgraph of $\D$, and  $\D$ is again a subgraph of the transitive closure. 
Moreover, all \DAG s with the same reachability matrix have the same transitive reduction and the same transitive closure and, therefore, the same ancestors and descendants. 

The following is an immediate consequence of Theorem~\ref{le1:CMgraph}(b). 

\bco\label{rem:Bodot2}
The \mSEM\ $\bfx$ is also a \mSEM\ on the transitive closure of every \DAG\ with reachability matrix $\sgn(B)$. 
\eco

We use this corollary to obtain necessary and sufficient conditions on a \CM\ $B$ as in \eqref{maxlin}   to be the \CM\ of a \mSEM. 
In contrast to Theorem~\ref{fixedpointeq} and Corollary~\ref{equiDML}(a)  we do not require that $B$ belongs to a specific given \DAG. 

\bthe\label{Bodotsquared}
 Let $\bfx$  be a \ML\ model as in \eqref{maxlin} with \CM\  $B$ such that ${\sgn}(B)$ is the reachability matrix of some \DAG.
 Define 
  \begin{align*}
 A:={\rm diag}(b_{11},\ldots,b_{dd}),\quad  B_0: =\Big(\frac{b_{ij}}{b_{ii}} \Big)_{d\times d},  \quad\mbox{and}\quad
A^\tc_0: =B_0 - {\rm id}_{d\times d}, 
 \end{align*}
 where ${\rm id}_{d\times d}$ denotes the identity matrix.    Then $\bfx$ is a \mSEM\ if and only if  the following fixed point equation holds:
\begin{align}\label{fpe1}
B&=B\odot B_0,\quad \text{which is equivalent to}\quad  B= A\vee B\odot A_0^\tc.
\end{align}
\ethe

\bproof
Let  $\D^\tc$ be the transitive closure of a \DAG\ with node set $V=\{1,\ldots,d\}$ and reachability matrix  $\sgn(B)$. First we show that  $\bfx$ is a \mSEM\ if and only if the fixed point equation  $B=A \vee B \odot A^\tc_0$ holds. By Corollary~\ref{rem:Bodot2} $\bfx$ is a \mSEM\ if and only if it is a \mSEM\ on $\D^\tc$. 
Thus, by Theorem~\ref{fixedpointeq} it suffices  to show that  $A^\tc_0$ is equal to the weighted adjacency matrix $A_0=\Big(\frac{b_{ij}}{b_{ii}}\1_{\pa(j)}(i) \Big)_{d\times d}$ (the parents in $\D^\tc$) of $\D^\tc$.  
We denote by $\an(i)$ for $i\in V$ the  ancestors   of node $i$ in $\D^\tc$, and observe from the definition of $\D^\tc$ that $\an(i)=\pa(i)$ for all $i\in V$. 
Since $B_0$  is a weighted reachability matrix of $\D^\tc$, we obtain
\begin{align*}
A_0^\tc=B_0-{\rm id}_{d\times d}= \Big(\frac{b_{ij}}{b_{ii}}{\1}_{\an(j)}(i) \Big)_{d \times d}=\Big(\frac{b_{ij}}{b_{ii}}{\1}_{\pa(j)}(i) \Big)_{d \times d}.
\end{align*}
It remains to show that  $ B\odot B_0=A\vee B\odot A_0^\tc$. By the definition of the matrix product $\odot$ in \eqref{odot} the $ji$-th entry of $A\vee B \odot A^\tc_0$ is equal to
\begin{align*}
b_{ji}\1_{\{i\}}(j) \vee  \bigvee_{k=1 }^d b_{jk} \Big(\frac{b_{ki}}{b_{kk}} -\1_{\{i\}}(k)\Big) 
& =  b_{ji}\1_{\{i\}}(j) \vee \bigvee_{k=1\atop k\neq i }^d\frac{b_{jk}b_{ki}}{b_{kk}}\\
 &= b_{ji}\1_{\{i\}}(j) \vee b_{ji}\1_{V\setminus\{i\}}(j) \vee\bigvee_{k=1\atop k\neq i,j}^d \frac{b_{jk}b_{ki}}{b_{kk}}\\
 & = \bigvee_{k=1}^d \frac{b_{jk}b_{ki}}{b_{kk}},
\end{align*}
which is the $ji-$th entry of the matrix $B\odot B_0$. 
\eproof

A non-negative symmetric matrix is by Theorem~\ref{Bodotsquared} a \CM\ of a \mSEM\ if and only if it is  a weighted reachability matrix of a \DAG\ and satisfies \eqref{fpe1}. Assume that we have verified these properties for a matrix $B$.  In order to find now all \mSEM s  which have \CM\ $B$ we can first use \eqref{eq:dagDB} to derive the \MD\DAG\ $\D^B$ from $B$ and then Theorem~\ref{le1:CMgraph}(b) to find all \DAG s and weights as in \eqref{ml-sem} such that \eqref{bs-max} holds.

\section{Backward and forward information in a recursive max-linear model\label{s6}}

In this section we apply our previous results to investigate, which components in a given node set of $\D$ are relevant for maximal information on some other component. 

{We know already from Corollary~\ref{le:IneqXi} that  $X_i\le \frac{b_{ii}}{b_{il}}X_l$ for all $i\in V$ and $l\in \Des(i)$  so that for some node set $U\subseteq V$ and all $i\in V$,
\begin{align} \label{ineqU}
\bigvee_{j \in {\An}(i)\cap U} \frac{b_{ji}}{b_{jj}}X_j \le X_i \le \bigwedge_{l \in {\Des}(i)\cap U} \frac{b_{ii}}{b_{il}} X_l.
\end{align}}

{The values of the bounds in \eqref{ineqU} can often be found as the maximum and minimum over a smaller number of nodes. 
We illustrate this by the following example.}

\bexam[Continuation of  Examples~\ref{ex:maxlin} and~\ref{ex:ml}: bounds\label{ex:bounds}] \\
For $U=\{1,2\}$ and $i=4$ we find by \eqref{ineqU} the lower bound
\begin{align}\label{specbounds}
\frac{b_{14}}{b_{11}} X_1\vee \frac{b_{24}}{b_{22}} X_2 \le X_4. 
\end{align}
We discuss the lower bound in \eqref{specbounds} and distinguish between two cases. 

First assume that the path $[1\to 2 \to 4]$ is max-weighted, which is by Theorem~\ref{paththrough}(a) equivalent to $b_{14}=\frac{b_{12}b_{24}}{b_{22}}$. 
From Corollary~\ref{le:IneqXi} we obtain
\begin{align*}
\frac{b_{12}}{b_{11}}X_1 \le  X_2, 
\quad \text{equivalently} \quad \frac{b_{14}}{b_{11}}X_1 \le \frac{b_{24}}{b_{22}} X_2.
\end{align*}
Therefore, the lower bound of $X_4$ in \eqref{specbounds} is always
$\frac{b_{24}}{b_{22}} X_2$.

Now assume  that the path $[1\to 2 \to 4]$ is not max-weighted. 
Since this is the only path from $1$ to $4$ passing through node $2$, this is by Theorem~\ref{paththrough}(b)  equivalent to $b_{14}>  \frac{b_{12}b_{24}}{b_{22}}$. 
From the max-linear representation \eqref{ml-noise} of $X_1$ and $X_2$ we have $\frac{b_{24}}{b_{22}} X_2 < \frac{b_{14}}{b_{11}} X_1$ if and only if
\begin{align*}
\frac{b_{12}b_{24}}{b_{22}} Z_1 \vee b_{24} Z_2<b_{14} Z_1,\quad \text{equivalently} \quad 
b_{24} Z_2<b_{14}Z_1. 
\end{align*}
The event $\{b_{24} Z_2<b_{14} Z_1\}$ has positive probability, since $Z_1$ and $Z_2$ are independent with support $\R_+$, giving $\frac{b_{14}}{b_{11}} X_1$  as lower bound.
But also the event $\{\frac{b_{14}}{b_{11}} X_1\le \frac{b_{24}}{b_{22}} X_2\}$ has positive probability, giving the lower bound $\frac{b_{24}}{b_{22}} X_2$.
\eexam

A  node $j\in \An(i)\cap U$  is relevant for the lower bound in \eqref{ineqU} if  no max-weighted path from $j$ to $i$  passes through some other node in $U$. 
Observe that this includes the observation made in Example~\ref{ex:bounds}. 
The nodes in the upper bound of \eqref{ineqU} have a similar characterization. 
We present a formal definition of these particular ancestors and descendants,
characterize them below in Lemma~\ref{poss1}, and give  an example afterwards. 

\bde\label{lowhigh}
Let $U\subseteq V$ and $i\in V$. 
\begin{enumerate}
\item[(a)] 
We call a node $j\in{\An}(i)\cap U$  {\em lowest max-weighted ancestor of $i$ in $U$}, if no max-weighted  path from $j$ to $i$ passes through some node in $U\setminus\{j\}$. 
We denote the set of the lowest max-weighted ancestors of $i$ in $U$ by  ${\An}_{{\low}}^U(i)$. 
\item[(b)] 
We call a node $l\in{\Des}(i)\cap U$  {\em highest max-weighted descendant of $i$ in $U$}, if no max-weighted path from $i$ to $l$ passes through some node in $U\setminus\{l\}$.
 We denote the set of the highest max-weighted descendants of $i$ in $U$ by $\text{De}_{{\high}}^U(i)$. 
\end{enumerate}
 \vspace*{-0.5cm}
\ede 

For $i\in U$ we find that the only lowest max-weighted ancestor and the only highest max-weighted descendant of $i$ in $U$ is the node $i$ itself. 
 For $i\in U^c=V\setminus U$  a simple characterization of $\An_\low^U(i)$ and $\Des_\high^U(i)$  is given next; this allows us to identify these nodes via the \CM\ of  $\bfx$.

\ble \label{poss1} Let $U\subseteq V$ and $i\in V$.
\begin{enumerate}
\item[(a)] If $i\in U$, then $\An_\low^U(i)=\Des_\high^U(i)=\{i\}$.
\item[(b)] If  $i\in U^c$, then
\begin{align}
\An_\low^U(i)&=\Big \{j\in \an(i)\cap U: b_{ji}>
 \bigvee_{k \in \des(j)\cap U\cap\an(i)} \frac{b_{jk}b_{ki}}{b_{kk}}     \Big \} \label{anlow}\\
\Des_\high^U(i)&=\Big \{l\in \des(i)\cap U: b_{il}> 
\bigvee_{k \in \des(i)\cap U\cap\an(l)} \frac{b_{ik}b_{kl}}{b_{kk}}     \Big \}. \label{dehigh}
\end{align}
\end{enumerate}
\ele

\bproof
(a)\, {follows immediately from the definition.}\\[1mm]
(b)\, Since $i\in U^c$, we have by Definition~\ref{lowhigh}(a) that $\An_{\low}^U(i)\subseteq \an(i)\cap U$. For $j\in \an(i)\cap U$  
we know from Theorem~\ref{paththrough}(b)  that no max-weighted path from $j$ to $i$ passes through some node in $U\setminus\{j\}$ if and only if 
\begin{align*} 
b_{ji}> 
\bigvee_{k \in \des(j)\cap U\cap\an(i)} \frac{b_{jk}b_{ki}}{b_{kk}},
\end{align*}
where we have used for the equality  that $i\in U^c$. 
Similarly, we obtain \eqref{dehigh}.
\eproof 

\bexam\label{ex:bounds1}[Continuation of  Examples~\ref{ex:maxlin}, \ref{ex:ml}, \ref{ex:bounds}: ${\An}^U_{{\low}}(4)$] \\
In order to find the lowest max-weighted ancestors of node $4$ in $U=\{1,2\}$, first observe that 
the only max-weighted path $[2\to 4]$  from $2$ to $4$  does not pass through any node in $U\setminus\{2\}$. 
Therefore, we have  by Definition~\ref{lowhigh}(a) that 
$2 \in {\An}^U_{{\low}}(4)$.  
For node $1$ we consider -- as in  Example~\ref{ex:bounds} --  two cases and use \eqref{anlow}:
\begin{enumerate}
\item[(1)]
If $b_{14}=  \frac{b_{12} b_{24}}{b_{22}}$, then ${\An}^U_{{\low}}(4)=\{2\}$. 
\item[(2)]
If $b_{14}> \frac{b_{12} b_{24}}{b_{22}}$, then ${\An}^U_{{\low}}(4)=\{1,2\}$. 
 \end{enumerate}
Comparing this with Example~\ref{ex:bounds} shows that the {lower} bound of $X_4$  is indeed always realized by some lowest max-weighted ancestor of node 4 in~$U$. 
\eexam

{We prove  that the lower and upper bounds in \eqref{ineqU} are always realized by some  lowest max-weighted ancestor  and  highest max-weighted descendant in $U$, respectively.  
For the lower bound this is based on the fact that between all nodes in $\D$ and their ancestors in $U$ there is always a max-weighted path,  which contains a lowest max-weighted ancestor in $U$. }
For the upper bound we use the existence of a max-weighted path between all nodes and their descendants in $U$ that passes through some highest max-weighted descendant in $U$.  
Before we state the modified lower and upper bounds in Proposition~\ref{co:IneqXi1},  we provide a useful characterization for a path analysis, which includes these statements.

\ble \label{lem3}
Let $U\subseteq V$ and $i\in V$.
\begin{enumerate}
\item[(a)] 
$\D$ has a max-weighted path from $j$ to $i$ passing through some node in $U$ if and only if it has a max-weighted path from $j$ to $i$    passing through some node  in ${\An}^U_{{\low}}(i)$.
\item[(b)] 
$\D$ has a max-weighted path from $i$ to $l$ passing through  some node in $U$ if and only if it has a max-weighted path from $i$ to $l$    passing through some node  in  ${\Des}^U_{{\high}}(i)$.
\end{enumerate}
\ele

\bproof 
We only show (a), since (b) can be proved analogously.   
Assume that a max-weighted path from $j$ to $i$ passes through some node in ${\An}^U_{{\low}}(i)$. 
Since ${\An}^U_{{\low}}(i)\subseteq U$, there is obviously also a max-weighted path from $j$ to $i$ that passes through some node in $U$.    

{For the converse, we may assume that $i\in U^c$, since by Lemma~\ref{poss1}(a) ${\An}^U_{{\low}}(i)=\{i\}$ for $i\in U$ and hence every max-weighted path contains a node in ${\An}^U_{{\low}}(i)$.  
Among all max-weighted paths from $j$ to $i$ let $p$ be one with maximum number of nodes in $U$.} 
 Denote by  $k_1$  the lowest node on $p$  contained in $U$; i.e., the subpath of $p$  from $k_1$ to $i$ contains no other node of $U$.
Assume that $k_1\not\in \An^U_\low(i)$. 
Since $k_1\in U$ and $i\in U^c$, there is by Definition~\ref{lowhigh}(a) a max-weighted path $p_1$ from $k_1$ to $i$  that passes through some node $k_2\in U$ with $k_2\neq k_1$. 
Thus  by replacing in $p$ the subpath from $k_1$ to $i$ by $p_1$ we obtain by Remark~\ref{rem:mwp}(iii) a max-weighted path from $j$ to $i$ containing more nodes in $U$ than $p$. 
This is however a contradiction.
Hence,  $k_1\in {\An}_{{\low}}^U(i)$, and  $p$ is a max-weighted path from $j$ to $i$ that passes through some node in ${\An}_{{\low}}^U(i)$. 
\eproof

\bpr \label{co:IneqXi1}
Let $U\subseteq V$ and $i\in V$. Then
\begin{align}\label{ineqU1}
\bigvee_{j \in {\An}(i)\cap U} \frac{b_{ji}}{b_{jj}}X_j=\bigvee_{j \in {\An}^U_{{\low}}(i)} \frac{b_{ji}}{b_{jj}} X_j \quad \text{and}\quad  \bigwedge_{l \in {\Des}(i)\cap U} \frac{b_{ii}}{b_{il}} X_l= \bigwedge_{l \in {\Des}^U_{{\high}}(i)} \frac{b_{ii}}{b_{il}} X_l .
\end{align}
\epr

\bproof 
Note from Definition~\ref{lowhigh}(a) that ${\An}^U_{{\low}}(i)\subseteq \An(i)\cap U$. To show the first equality  take some $k\in ({\An}(i)\cap U)\setminus {\An}^U_{{\low}}(i)$. Observe from   Lemma~\ref{poss1}(a) that $k\neq i$ and, hence, $k\in \an(i)\cap U$.  By  Lemma~\ref{lem3}(a) there must be a max-weighted path from $k$ to $i$, which passes through  some node  $j\in {\An}_{{\low}}^U(i)$.  
By \eqref{paththrough2} and  Corollary~\ref{le:IneqXi}, we obtain
\begin{align}\label{eq:IneqXi1}
\frac{b_{ki}}{b_{kk}}X_k =  \frac{b_{kj}b_{ji}}{b_{kk}b_{jj}} X_k \le \frac{b_{ji}}{b_{jj}}X_j.
\end{align}
Since for all $k \in  ({\An}(i)\cap U)\setminus {\An}^U_{{\low}}(i)$ there exists some $j \in {\An}_{{\low}}^U(i)$ such that \eqref{eq:IneqXi1} holds, the first equality of \eqref{ineqU1} follows. 
The second equality may be verified analogously. 
\eproof

So far, for every component  of $\bfx$, we have identified  a lower and upper bound in terms of the components of $\bfx_U=(X_l, l\in U)$. 
However, we cannot say anything about the quality of the bounds. 
For instance, we do not know in which situation a component attains one of the bounds. 
We clarify this by writing all components of $\bfx$ 
as max-linear functions of $\bfx_U$ and certain noise variables. There are many such representations, since we can always include non-relevant ancestral components with appropriate \ML\ coefficients as we know from Theorem~\ref{le1:CMgraph}(b).
To find the relevant components  of $\bfx_U$ and noise variables we focus on those with the minimum number of components of $\bfx_U$ and the minimum number of noise variables. 
For $i\in V$ we denote by $\an^U_\nmw(i)$ {the set of all $j\in\an(i)$ 
such that no max-weighted path from $j$ to  $i$ passes through some node in $U$. }
By Theorem~\ref{paththrough}(b) we have 
\begin{align}\label{nmw}
\an_{\nmw}^U(i)=\Big\{j\in \an(i): b_{ji}> \bigvee_{k\in \Des(j)\cap U\cap \An(i)} \frac{b_{jk}b_{ki}}{b_{kk}}\Big\}. 
\end{align}
{Since  $j\in \an(i)\setminus\an^U_\nmw(i)$ if and only if there is a max-weighted path from $j$ to $i$ passing through some node in $U$, we have by Theorem~\ref{paththrough}(a)
\begin{align}\label{mw}
\an(i)\setminus\an_{\nmw}^U(i)=\Big\{j\in \an(i): b_{ji}= \bigvee_{k\in \Des(j)\cap U\cap \An(i)} \frac{b_{jk}b_{ki}}{b_{kk}}\Big\}. 
\end{align}}

\bthe \label{repVB} 
Let $\bfx$ be a \mSEM\ with \DAG\ $\D$ and \CM\ $B$, and let  $U\subseteq V$. 
Let ${\An}^U_{{\low}}(i)$ be the lowest max-weighted ancestors of node $i$  in $U$ as in Definition~\ref{lowhigh}(a),  and define $\An^{U}_{{\nmw}}(i):=({\an}_{{\nmw}}^U(i)\cup\{i\})\cap U^c$.  
Then for every $i\in V$,
 \begin{align}\label{eqrep}
 X_i &=  \bigvee_{k \in {\An}^U_{{\low}}(i)} \frac{b_{ki}}{b_{kk}} X_k \vee   \bigvee_{j \in \An^{U}_{{\nmw}}(i)} b_{ji}Z_j. 
\end{align} 
This representation of $X_i$ as a max-linear function of $\bfx_U$ and noise variables involves  the minimum number of components of $\bfx_U$ and the minimum number of   noise variables. 
 \ethe

 \bproof
 We distinguish between nodes $i\in U$ and $i\in U^c$. 
 For $i\in U$ we know from  Lemma~\ref{poss1}(a) that $\An^U_{{\low}}(i)=\{i\}$. 
 Furthermore, we  have $\An_{\nmw}^U(i)=\emptyset$, since $i\in U$ and every path, hence also every max-weighted  path, from some $j\in\an(i)$ to $i$ passes through some node in $U$, namely $i$ itself. 
 Thus we obtain \eqref{eqrep}. The second statement is obvious. 

Now assume that $i\in U^c$, and note that in this case $\An_{\nmw}^U(i)=\an_{\nmw}^U(i)\cup\{i\}$.
Applying the first equality in \eqref{ineqU1} and \eqref{ml-noise} as well as \eqref{lem14} in a second step to interchange the first two maximum operators, we have 
\begin{align} 
 \bigvee_{k \in {\An}^U_{{\low}}(i)} \frac{b_{ki}}{b_{kk}} X_k &= \bigvee_{k \in{\an}(i)\cap U} \frac{b_{ki}}{b_{kk}} \big(\bigvee_{j \in \An(k)} b_{jk}Z_j\big)=  \bigvee_{j\in {\an}(i)}  \bigvee_{k \in {\Des}(j)\cap  {\an}(i)\cap U} \frac{b_{jk}b_{ki}}{b_{kk}} Z_j  \label{eqrep4}. 
\end{align}
We split up the set ${\an}(i)$ into ${\an}_{{\nmw}}^U(i)$  and ${\an}(i)\setminus {\an}_{{\nmw}}^U(i)$ as well as the set $\An^{U}_{{\nmw}}(i)$ into ${\an}_{{\nmw}}^U(i)$ and $\{i\}$ to obtain that the right-hand side of \eqref{eqrep} is equal to
\begin{align*}
&\bigvee_{j\in \an(i)\setminus {\an}_{{\nmw}}^U(i)}  \bigvee_{k \in {\Des}(j)\cap {\an}(i)\cap U} \frac{b_{jk}b_{ki}}{b_{kk}}  Z_j \\
&\quad\quad\quad\quad \vee \bigvee_{j\in {\an}^U_{{\nmw}}(i)} \big( \bigvee_{k \in {\Des}(j)\cap {\an(i)}\cap U} \frac{b_{jk}b_{ki}}{b_{kk}} \vee b_{ji}\big) Z_j \vee  b_{ii}Z_i.
\end{align*}
Noting that $i\in U^c$ when using   \eqref{mw} and \eqref{nmw} yields for the right-hand side of \eqref{eqrep} 
\begin{align*}
 \bigvee_{j\in \an(i)\setminus {\an}_{{\nmw}}^U(i)}  b_{ji} Z_j \vee   \bigvee_{j \in \an^{U}_{{\nmw}}(i)} b_{ji}Z_j \vee b_{ii}Z_i= \bigvee_{j \in \An(i)} b_{ji}Z_j=X_i.
\end{align*}
In order to verify that for $i\in U^c$  \eqref{eqrep} is the representation of $X_i$ with the minimum number of components of $\bfx_U$ and the minimum number of noise variables, we prove that each term on the right-hand side of \eqref{eqrep} has to appear, since  otherwise some noise variable $Z_j$ in representation \eqref{ml-noise} would have a weight strictly less than $b_{ji}$. 
{We compare the noise variables of the right-hand sides of \eqref{eqrep} and  \eqref{eqrep4}.}
Since    $b_{ii}Z_i$ does not appear in \eqref{eqrep4}, it has to to appear in \eqref{eqrep}. 
For $j\in\an^U_{\nmw}(i)$ it follows from \eqref{nmw}  that if $Z_j$ appears in \eqref{eqrep4}, then with a coefficient strictly less than $b_{ji}$. 
The maximum over $\An^{U}_{{\nmw}}(i)$ must therefore appear in \eqref{eqrep}.    
{Definition~\ref{lowhigh}(a) implies that no max-weighted path from $j\in\An_{\low}^U(i)$ to $i$ passes through some node in $\des(j)\cap U\cap \An(i)$.}
Thus observe from \eqref{eqrep4} and \eqref{paththrough3}  that only the term $\frac{b_{ji}}{b_{jj}}X_j$ provides $Z_j$ with the weight $b_{ji}$ on the right-hand side of  \eqref{eqrep} and the term $\frac{b_{ji}}{b_{jj}}X_j$  has to appear on the right-hand side of \eqref{eqrep}. 
 \eproof


{We use Theorem~\ref{repVB} to obtain for every component $X_i$ a minimal representation in terms of the components of $\bfx_{\pa(i)}$ and independent noise variables.}

\bco\label{cor:minrep}
Let  $\D^B$ be the  \MD \DAG\   of $\bfx$ as in Definition~\ref{DB} with parents $\pa^B(i)$ of node $i$ in $\D^B$.
Then  for all $i\in V$ we have 
${\An}^{{\pa}(i)}_{{\low}}(i)={\pa}^B(i)$ and 
\begin{align}\label{minrep}
X_i&= \bigvee_{k \in {\pa}^B(i)} \frac{b_{ki}}{b_{kk}}X_k \vee b_{ii}Z_i= \bigvee_{k \in {\pa}^B(i)} c_{ki}X_k \vee c_{ii} Z_i. 
\end{align}
\eco

\bproof
Recall from 
\eqref{DB1} that
\begin{align*}
\pa^B(i)=\Big\{ k\in \pa(i): b_{ki}> \bigvee_{l\des(k)\cap \pa(i)} \frac{b_{kl}b_{li}}{b_{ll}} \Big\},
\end{align*}
and observe from this and \eqref{anlow} that ${\An}^{{\pa}(i)}_{{\low}}(i)={\pa}^B(i)$.   
Since every path from $j\in\an(i)$ to $i$  passes through some node in $\pa(i)$, there  is always a max-weighted path from $j$ to $i$ containing some node of $\pa(i)$.
Hence, $\An_{{\nmw}}^{{\pa}(i)}(i)=(\an_\nmw^{\pa(i)}(i)\cup\{i\})\cap(\pa(i))^c=\{i\}$.
Thus  we obtain by Theorem~\ref{repVB}   the first equality  in \eqref{minrep}. 
For the second, recall from Theorem~\ref{le1:CMgraph}(a) that $c_{ki}=\frac{b_{ki}}{b_{kk}}$ for $k \in {\pa}^{B}(i)$. 
\eproof


\brem 
Representation \eqref{minrep}  complements Theorem~\ref{le1:CMgraph}(a); we find again that the  \MD \DAG\  $\D^B$ yields the
minimal representation of  $\bfx$ as a \mSEM. 
\erem

The following example  illustrates Theorem~\ref{repVB}. 

\bexam\label{ex:bounds2}[Continuation of  Examples~\ref{ex:maxlin},~\ref{ex:ml},~\ref{ex:bounds}, and~\ref{ex:bounds1}:  minimal representation of $X_4$ by $\bfx_U$] \\
We consider again $U=\{1,2\}$ and $i=4$.  
Obviously, there are max-weighted paths from $1$ and $2$ to $4$ passing through some node in $U=\{1,2\}$. 
Hence, $1,2 \in \an(4)\setminus \an_\nmw^U(4)$.   
Since no max-weighted path from $3$ to $4$ passes through $1$ or $2$, we have $\An^U_{{\nmw}}(4)=(\an^U_{{\nmw}}(4)\cup\{4\})\cap U^c=\{3,4\}$. 
In Example~\ref{ex:bounds1} we have already determined the set $\An^U_{\low}(4)$ depending on the \ML\ coefficients.   
Thus we distinguish again between  two cases:
\begin{enumerate}
\item[(1)]
If $b_{14}=  \frac{b_{12} b_{24}}{b_{22}}$, then $X_4=\frac{b_{24}}{b_{22}}X_2\vee b_{34}Z_3\vee b_{44}Z_4$.\\
We want to remark that the conditional independence properties of $\bfx$ are reflected in this representation:  from Example~\ref{ex:ml} we know  that $X_1\upmodels X_4\mid X_2$. So it is obvious that  $X_1$ does not appear in the minimal representation of $X_4$ as max-linear function of $X_1$ and $X_2$.
\item[(2)]
If $b_{14}> \frac{b_{12} b_{24}}{b_{22}}$, then  $X_4=\frac{b_{14}}{b_{11}}X_1\vee \frac{b_{24}}{b_{22}}X_2\vee b_{34}Z_3\vee b_{44}Z_4$.\\
In particular, $ \frac{b_{14}}{b_{11}}X_1>\frac{b_{24}}{b_{22}}X_2$ is possible with positive probability; in (1) this is not possible  (see Example~\ref{ex:bounds}). 
\end{enumerate}
In both representations all random variables have to appear, but no other ones are needed.
Hence, we have indeed derived the minimal representation of $X_4$ in terms of $X_1$ and $X_2$. 

For $U=\{2\}$ and $i=4$ we have ${\An}^U_{{\low}}(4)=\{2\}$. 
Similarly as above we obtain that  $2\in \an(4)\setminus \an_\nmw^U(4)$ and  $3,4 \in  \An^U_{{\nmw}}(4)$.  
It remains to discuss node $1$, which gives rise to the same two cases as above:
\begin{enumerate}
\item[(1)]
If the path $[1\to 2\to 4]$ is max-weighted, then
$X_4=\frac{b_{24}}{b_{22}}X_2\vee b_{34}Z_3\vee b_{44}Z_4.$
\item[(2)]
If the path $[1\to 2\to 4]$ is not max-weighted, then
$X_4=\frac{b_{24}}{b_{22}}X_2\vee b_{14}Z_1\vee b_{34}Z_3\vee b_{44}Z_4.$
\end{enumerate} 
{Such minimal representations become relevant, when $\bfx$ is partially observed.
If, for example, $X_2$ is observed, then the prediction problem of $X_4$ can be solved by the observations of $X_2$ and by conditional simulation of  the relevant noise variables; see \cite{Wang2011}.
In case (1) we need to simulate $Z_3, Z_4$, whereas in case (2) $Z_1,Z_3, Z_4$ are needed. 
We will discuss such prediction problems in a follow-up paper.}
\eexam

\appendix

\section{Auxiliary lemma\label{aux}} 

\ble\label{lem1}
Let $\D=(V,E)$ be a \DAG\ and $U\subseteq V$. 
For non-negative functions $a(i,j,k)$ for  $i,j,k \in V$  we have for all $i\in V$,
\begin{align}
\bigvee_{k\in{\pa}(i)} \bigvee_{j\in\an(k)} a(i,j,k)  &= \bigvee_{j\in{\an}(i)} \bigvee_{k\in{\des}(j)\cap {\pa}(i)}  a(i,j,k) \label{lem11}\\
\bigvee_{k \in {\an}(i)\cap U }\bigvee_{j \in \An(k)}  a(i,j,k) &=   \bigvee_{j \in {\an}(i)} \bigvee_{k \in {\Des}(j)\cap {\an}(i)\cap U  }    a(i,j,k). \label{lem14}
\end{align}
\ele

\bproof 
Since we take maxima,  we only have to prove that each combination of nodes $(k,j)$ on the left-hand side appears also on the right-hand side and vice versa. 
In order to prove \eqref{lem11}, 
it suffices to show that
\begin{align*}
k \in {\pa}(i) \text{ and } j \in \an(k) \quad \text{if and only if} \quad j \in {\an}(i) \text{ and } k \in {\des}(j)\cap {\pa}(i).
\end{align*}
By observing that $\an(\pa(i))\subseteq \an(i)$ and $j \in \an(k)$ if and only if $k \in {\des}(j)$ this equivalence is obvious. 
Eq. \eqref{lem14} is  proved in the same way. 
\eproof

\section*{Acknowledgements}
We thank Steffen Lauritzen for fruitful discussions and his constructive suggestions, which improved our manuscript. 
Nadine Gissibl had the pleasure of spending two months at the Seminar for Statistics of the ETH Zurich. She wants to thank all colleagues there for a very pleasant time. She also gratefully  acknowledges support from the TUM Graduate School's {\em International School of Applied Mathematics}.

\bibliographystyle{abbrvnat}
\bibliography{Graphs}

\end{document}